\let\ifiso\iffalse
\let\ifdump\iffalse
\catcode`\{=1 
\catcode`\}=2 
\catcode`\$=3 
\catcode`\&=4 
\catcode`\#=6 
\catcode`\^=7 \catcode`\^^K=7 
\catcode`\_=8 \catcode`\^^A=8 
\catcode`\^^I=10 
\chardef\active=13 \catcode`\~=\active 
\catcode`\^^L=\active \outer\def^^L{\par} 
\catcode`\@=11
\let\e@\expandafter
\def\s@tt#1#2{\e@\e@\e@\let\e@#1#2}
\def\letifundefined#1#2{\def\aaa{#1}\def\bbb{#2}%
\ifx#1\undefined\s@tt\aaa\bbb\fi}
\letifundefined\ifdump\iftrue
\letifundefined\ifmac\iffalse
\letifundefined\ifiso\iftrue
\catcode`\@=12

\ifmac
\input option_keys
\else
\ifdump
\scrollmode
\ifiso
\input isoplain
\def\v#1{\char\if c#1"A3\else\if C#1"83\else\if s#1"92\else\if
S#1"B2\else\if d#1"A4\else\if D#1"84\else\if E#1"85\else\if
#1e"A5\else\if Z#1"9A\else\if z#1"BA\else\if N#1"8C\else\if n#1"AC
\fi\fi\fi\fi\fi\fi\fi\fi\fi\fi\fi\fi}
\else
\input plain
\fi
\fi
\fi
\catcode`\@=11 
\ifiso
\else

\let\ifmac\iffalse
\toksdef\t@=1
\let\e@\expandafter
\def\oms@nn{\iftrue}
\def\omf@lsk{\iffalse}
\def\ifdefined#1{\e@\ifx\csname
#1\endcsname\relax\e@\omf@lsk\else\e@\oms@nn\fi}
\def\newl@nguage{\alloc@9\language\chardef\@cclvi}

\def\makelanguage#1{\ifdefined{l@ng#1}
\message{The language #1 is already defined.}%
\else
\e@\newl@nguage\csname l@ng#1\endcsname
\e@\language\csname l@ng#1\endcsname
\language0
\fi}

\def\makel@nguage#1{\ifdefined{l@ng#1}
\message{The language #1 is already defined.}%
\else
\message{Making language #1.}
\e@\newl@nguage\csname l@ng#1\endcsname
\e@\language\csname l@ng#1\endcsname
\csname input\endcsname #1hyphen
\language0
\fi}

\def\deflanguage#1#2{\ifdefined{l@ng#1}%
\e@\def\csname #1\endcsname{\e@\language\csname l@ng#1\endcsname#2}\else
\errmessage{The language #1 is not defined.}\fi}

\def\l@ngdef#1#2#3{\def#1{\csname \string#1@\the\language\endcsname}%
\edef\tmp@{\e@\the\csname l@ng#2\endcsname}%
\e@\def\csname \string#1@\tmp@\endcsname{#3}}

\def\langdef#1{\t@{#1}\def\l@ngnext\in##1=##2##3{\if##3.\let\l@ngnext\relax\fi
\e@\l@ngdef\the\t@{##1}{##2}\l@ngnext}\l@ngnext}

\let\m@kelang\makelanguage
\ifdump
\ifiso
\ifmac
\else
\let\m@kelang\makel@nguage
\fi
\fi
\fi

\m@kelang{english}
\deflanguage{english}{\lefthyphenmin=2\righthyphenmin=3\nonfrenchspacing}

\m@kelang{swedish}
\deflanguage{swedish}{\lefthyphenmin=2\righthyphenmin=2\frenchspacing}

\m@kelang{icelandic}
\deflanguage{icelandic}{\lefthyphenmin=2\righthyphenmin=2\nonfrenchspacing}

\m@kelang{french}
\deflanguage{french}{\lefthyphenmin=2\righthyphenmin=2\frenchspacing}

\langdef\today\in english=\Today,\in swedish=\idag,%
\in french=\aujourdhui,\in icelandic=\iidag.

\def\Today{\ifcase\month\or
 Jan\or Febr\or  Mar\or  Apr\or May\or Jun\or  Jul\or
 Aug\or  Sep\or  Oct\or Nov\or  Dec\or\fi
 \space\number\day, \number\year}

\def\idag{\number\day\space\ifcase\month\or
 januari\or februari\or  mars\or  april\or maj\or juni\or  juli\or
 augusti\or  september\or  oktober\or november\or  december\or\fi,
\number\year}

\def\aujourdhui{le \number\day\space\ifcase\month\or
 janvier\or f{\`e}vrier\or  mars\or  avril\or mai\or juin\or  juillet\or
 aôut\or  septembre\or  octobre\or novembre\or  d{\`e}cembre\or\fi,
\number\year}

\def\iidag{\number\day.\space\ifcase\month\or
janúar\or febrúar\or mars\or apríl\or maí\or júní\or
júlí\or ágúst\or september\or október\or nóvember\or
desember\fi,\space\number\year}
\ifmac\else
\ifiso
\let\ídag\iidag
\fi
\fi

\english

\fi
\def\inputonce#1 {\e@\ifx\csname @np#1\endcsname\relax\input#1
\e@\let\csname @np#1\endcsname1\else\fi}
\def\e@t#1{}
\def\e@tif#1#2{\if#1#2\e@\e@t\fi}
\def\makesplit#1#2{\def#1##1{%
\def\@split####1#2####2#2{\def\firstsplit{####1}\def\secondsplit
{####2}\e@tif#2}%
\@split##1#2#2\relax}}
\def\eat#1{}

\def\@strue{\let\tr@th\iftrue}
\def\@sfalse{\let\tr@th\iffalse}
\def\boolnot#1{#1\@sfalse\else\@strue\fi\tr@th}

\def\letctrl#1\to#2 {\e@\let\csname #1\endcsname#2}

\def\setbool#1{#1\@struth\else\@sfalse\fi}

\def\ifdef#1{\edef\s@tt{#1}%
\boolnot{\e@\e@\e@\ifx\e@\csname\s@tt\endcsname\relax}}

\def\ifempty#1{\@ifempty #1\@emptymarkA\@emptymarkB}%
\def\@ifempty#1#2\@emptymarkB{\ifx #1\@emptymarkA}%
\def\@emptymarkA{\@emptymarkA}

\def\omtom#1{\ifx\relax#1\relax}

\def\ifstreq#1#2{\edef\s@tt{#1}\edef\t@tt{#2}\ifx\s@tt\t@tt}

\scrollmode

\let\e@\expandafter

\newif\ifch@pq 
\newif\ifrefdef 
\newif\ifnoreferr
\newif\if@nd
\newif\iffinal
\newif\ifsvensk
\newif\iffransk
\newif\ifchap
\newif\ifjour@
\newif\ifbook@
\newif\ifinbook@
\newif\ifpages@
\newif\ifpage@
\newif\ifsln@
\newif\ifprep@
\newif\ifspec@
\newif\iflect@
\newif\ifbibliography
\newif\iffriendly
\newif\ifst@rted 

\newdimen\las
\newcount\chapno
\newcount\sekt
\newcount\lop
\newcount\ovnnr

\def\newtoks@{\alloc@5\toks\toksdef\@cclvi} 
\long\def\rightaddtoks#1#2{#1\e@{\the#1#2}}
\long\def\leftaddtoks#1#2{\t@{#2}#1\e@\e@\e@{\e@\the\e@\t@\the#1}}
\def\tokconcat#1=#2&#3{#1\e@\e@\e@{\e@\the\e@#2\the#3}}

\def\regs@t#1#2{\e@\let\e@#1\csname#2\endcsname}

\def\mktokreg#1{\e@\newtoks@\csname #1\endcsname}
\def\tokrightappenditem#1\to#2 {\e@\rightaddtoks\csname #2\endcsname{#1}}
\def\tokleftappenditem#1\to#2 {\e@\leftaddtoks\csname #2\endcsname{#1}}
\def\reftokreg#1{\e@\the\csname #1\endcsname}

\def\tokregconcat#1=#2&#3 {%
\regs@t\r@tt{#1}%
\regs@t\s@tt{#2}%
\regs@t\t@tt{#3}%
\r@tt\e@\e@\e@{\e@\the\e@\s@tt\the\t@tt}}

\toksdef\t@=1  
\toksdef\t@@=2 
\newtoks\byes
\newtoks\presentpath
\newcount\x@ant
\newtoks\x@note
\newtoks\x@sid
\newtoks\x@bib

\long\def\left@ppenditem#1#2{\t@={\\{#2}}\t@@=\e@{#1}%
\edef#1{\the\t@\the\t@@}}
\long\def\right@ppenditem#1#2{\t@={\\{#2}}\t@@=\e@{#1}%
\edef#1{\the\t@@\the\t@}}
\def\mkreg#1{\e@\def\csname#1\endcsname{}}
\long\def\leftappenditem#1\to#2 {\e@\left@ppenditem\csname#2\endcsname{#1}}
\long\def\leftappenditem#1\to#2 {\e@\right@ppenditem\csname#2\endcsname{#1}}
\def\refreg#1{\csname #1\endcsname}

\def\setr@g#1\to#2{\e@\let\csname #1\endcsname#2}

\outer\def\Bye{\supereject\@nd}

\outer\def\Beginsection#1\par{\advance\sekt by 1\lop=0
\vskip \z@ plus.3\vsize \penalty -250 \vskip \z@
plus-.3\vsize \bigskip \vskip \parskip \message {\number\sekt. #1}\leftline
{\bf\number\sekt.  #1}\nobreak \smallskip \noindent}

\outer\long\def\Proclaim#1. {\medbreak
 \noindent{\bf#1.\enspace}\bgroup\sl
\def\pro{\egroup\par
 \ifdim\lastskip<\medskipamount \removelastskip\penalty55\medskip\fi
 {\bf \sevenrm\ifsvensk BEVIS\else PROOF\fi}: }}

\outer\long\def\Proclame#1. {\medbreak
 \noindent{\bf#1.\enspace}\bgroup\sl
\def\pro{\egroup\par
 \ifdim\lastskip<\medskipamount \removelastskip\penalty55\medskip\fi}}

\def\eqalignno#1{\global\ch@pqtrue
\displ@y \tabskip \centering \halign to\displaywidth
{\hfil $\displaystyle {##}$\tabskip \z@skip &$\displaystyle {{}##}$\hfil
\tabskip\centering &\t@={##}%
\omtom{\the\t@}%
\else
 \global\advance\lop by 1%
\global\e@\edef\csname l@\e@ p\the\t@\endcsname
{\ifchap\the\chapno @\fi\iffinal\sektname.\number\lop\else\the\t@\fi}
\llap {(\iffinal\sektname.\number\lop\else\the\t@\fi)}
\fi%
\tabskip \z@skip \crcr #1\crcr }}

\def\leqalignno#1{\global\ch@pqtrue
\displ@y \tabskip \centering \halign to\displaywidth
{\hfil $\displaystyle {##}$\tabskip \z@skip &$\displaystyle {{}
##}$\hfil \tabskip \centering &\t@={##}\omtom{\the\t@}\else
\global\advance\lop by 1
\global\e@\edef\csname l@\e@ p\the\t@\endcsname
{\ifchap\the\chapno @\fi\iffinal\sektname.\number\lop\else
\the\t@\fi}
\kern -\displaywidth \rlap
{(\iffinal\sektname.\number\lop\else\the\t@\fi)}\fi
\tabskip \displaywidth \crcr #1\crcr}}

\font\dummyft@=dummy
\def\fontlist@{\\{\tenrm}\\{\sevenrm}\\{\fiverm}\\{\teni}\\{\seveni}%
 \\{\fivei}\\{\tensy}\\{\sevensy}\\{\fivesy}\\{\tenex}\\{\tenbf}\\{\sevenbf}%
 \\{\fivebf}\\{\tensl}\\{\tenit}\\{\frak}\\{\tentt}%
 }
\def\font@#1=#2 {\rightappend@#1\to\fontlist@\font#1=#2 }
\def\dodummy@{{\def\\##1{\let##1\dummyft@}\fontlist@}}
\newif\ifsyntax@
\newcount\countxviii@
\def\nopages@{\output={\setbox\z@\box255 \deadcycles\z@}\newtoks@\output}
\def\syntax{\hfuzz3000pt \vfuzz3000pt\tolerance\@M
\syntax@true\dodummy@\countxviii@\count18
 \loop\ifnum\countxviii@>\z@\textfont\countxviii@=\dummyft@
 \scriptfont\countxviii@=\dummyft@\scriptscriptfont\countxviii@=\dummyft@
 \advance\countxviii@\m@ne\repeat
 \dummyft@\tracinglostchars\z@\nopages@\frenchspacing\hbadness\@M}

\font\ninerm=cmr9

\font\frak=eufm10 
\font\smc=cmcsc10 
\def\skrivstart#1{%
\e@\let\e@\t@mp\csname newwrite\endcsname
\e@\t@mp\csname stre@m#1\endcsname
\e@\immediate\e@\openout\csname stre@m#1\endcsname=\rjobname.#1 \par}

\def\skriv#1{\e@\immediate\e@\write\csname stre@m#1\endcsname}

\def\skrivslut#1{\e@\immediate\e@\closeout\csname stre@m#1\endcsname}

\def\startout#1#2{\skrivstart{#1}\e@\let\csname #1out\endcsname1
\unitaction{#1}{\t@={}\def\t@steg####1{#2####1{#1}}%
\unittotoken{#1}\t@steg}%
\addafterunit{document}{\skrivslut{#1}}}

\def\writeandput#1#2{{\newlinechar=`\ %
\skriv{#2}{\noexpand\ovnnr\the\ovnnr\noexpand\begin{#2}\the#1}}\the#1}

\def\writeonly#1#2{{\newlinechar=`\ %
\skriv{#2}{\noexpand\ovnnr\the\ovnnr\noexpand\begin{#2}\the#1}}\egroup}

\def\mdef#1#2{\def#1{\relax\ifmmode {#2}\else$#2$\fi}} 

\def\acce{\relax\ifmmode\acute e\else {\'e}\fi}

\def\rum{\nonscript\!}
\let\@a=\in
\def\in{\rum\@a\rum}
\mdef\cL{{\cal L}}

\mdef\infinity{\infty}
\def\hatt#1{\relax\ifmmode \hat#1\else \^#1\fi}
\def\tjeck#1{\relax\ifmmode \check#1\else\v#1\fi}
\def\prickar{\relax\ifmmode \ldots\else \dots\fi}

\mdef{\ovl#1}{\overline {#1}}
\mdef{\undl#1}{\underline {#1}}

\mdef{\ko#1}{{\cal O}_{#1}}

\def\icmslt@st#1#2#3#4#5#6@{\if#1d\if#2c\if#3s\if#4l}
\def\ncmslt@st#1#2#3#4#5#6@{\if#1c\if#2m\if#3s\if#4l}
\ifiso
\let\cmslt@st\icmslt@st
\else
\let\cmslt@st\ncmslt@st
\fi
\def\deffont{\let\t@mp\sl\e@\cmslt@st \fontname
      \font\ \ \ \ @\let\t@mp\rm\fi\fi\fi\fi\t@mp}

\def\definition#1{{\deffont #1}}

{\def\tag#1#2#3@{#1#2}
\def\\#1 {\if1#1\else\global\e@\mdef\csname
g\tag#1@\endcsname{\csname#1\endcsname}\e@
\\\fi}
\\alpha beta gamma Gamma delta Delta zeta eta theta Theta
iota kappa lambda Lambda mu nu rho sigma Sigma pi Pi
xi Xi chi  tau phi Phi omega Omega epsilon psi Psi upsilon Upsilon 1 }
{\def\tag var#1#2@{#1}
\def\\#1 {\if1#1\else\global\e@\mdef\csname
gva\tag#1@\endcsname{\csname#1\endcsname}\e@
\\\fi}
\\varrho varphi varpi vartheta varsigma 1 }
\mdef\gep\varepsilon

\mdef\Sp{{\mathop{\bf Spec\hskip1.5pt}\nolimits}}
\mdef\Spf{{\mathop{\bf Specf\hskip2pt}\nolimits}}
\mdef\Proj{{\mathop {\bf Proj\hskip1.5pt}\nolimits}}
\mdef\l{\ell}

\mdef{\pow#1}{{\bf k}\bigl[[#1]\bigr]}

\def\=#1{\relax\ifmmode\mathrel{\mathop=^{#1}}\else{\accent 22 #1}\fi}

\def\Dsum{\bigoplus}  

\let\co=\colon

\mdef\mul{{\bf G}_m}
\mdef\add{{\bf G}_a}

\mdef{\alt#1}{\alpha_{#1}}
\mdef{\Alt#1}{\alpha_{#1}}

\mdef{\Hom#1}{\underline{Hom}_{{\cal O}_{#1}}}
\mdef{\set#1#2}{\{#1 : #2\}}

\let\@d=\ni
\def\ni{\rum\@d\rum}

\def\riso{\mathrel{\hskip7pt\raise-2.5pt\hbox{$\widetilde{\phantom{xx}}$}
\kern-16pt\longrightarrow}}
\def\liso{\mathrel{\hskip7pt\raise-2.5pt\hbox{$\widetilde{\phantom{xx}}$}
\kern-17pt\longleftarrow}}

\def\subsetneq{\mathrel{\kern6pt\raise-2pt\hbox{$\scriptscriptstyle\not$}
\kern-5pt\subseteq}}
\def\supsetneq{\mathrel{\kern6pt\raise-2pt\hbox{$\scriptscriptstyle\not$}
\kern-3pt\supseteq}}

\def\disjunion{\mathbin{\vbox{\hbox{\kern1.8pt.}\kern-5.4pt\hbox{$\cup$}}}}

\def\fra#1{\relax\ifmmode\hbox{\frak #1}\else{\frak #1}\fi}

\mdef{\hod#1#2#3#4}{\if#30 H^{#2}({#1},{\cal O}_{#1}) \else
 H^{#2}(#1,\Omega^{#3}\if\relax#4\relax_{#1}\else _{#1/#4}\fi)\fi}

\mdef{\coh#1#2#3}{H^{#1}(#2,#3)}

\mdef{\Cal#1}{{\cal #1}}

\def\can@#1{\ifx#1\st@pit\else \e@\def\csname
#1\endcsname{{\bf #1}}\e@\can@\fi}
\can@  QPZNARCFk\st@pit
\mdef\bk{\bar{\bf k}} \mdef\bQ{\bar{\bf Q}}

\def\candef#1{\can@ #1\st@pit}

\let\endfrac\relax
\def\contfrac{\t@={}\t@@={}\contfr@c}
\def\contfr@c#1{%
\if#1\endfrac
 \the\t@
\else
 \t@=\e@{\the\t@\egroup}%
 \the\t@@\t@@={+}\bgroup\displaystyle\strut1\over\displaystyle #1\e@
 \contfr@c
\fi}

 \def\mapright#1{\setbox0=\hbox spread 12pt{\hfil$\scriptstyle#1$\hfil}
  \mathop{\hbox to \wd0{\rightarrowfill}\kern-\wd0\raise6pt\box0}}

\def\mapleft#1{\setbox0=\hbox spread 12pt{\hfil$\scriptstyle#1$\hfil}
  \mathop{\hbox to \wd0{\leftarrowfill}\kern-\wd0\raise6pt\box0}}

\mdef{\pil#1#2#3}{{#1}\co{#2}\to{#3}}

\def\funk#1#2#3#4{\eqalign{#1&\to #2\cr#3&\mapsto#4\cr}}

\let\ign@re\empty
\let\in@nit\empty

\def\ifunit#1{\ifdef{ev@ry#1}}
\def\ifignore#1{\e@\ifx\csname ign@re#1\endcsname\ign@re}
\def\ifinunit#1{\e@\ifx\csname in@nit#1\endcsname\in@nit}
\def\ifpresentunit#1{\ifstreq{#1}\presentunit}

\def\f@rbidden#1#2{%
\errhelp{You cannot have unit #2 inside unit #1.}
\errmessage{Starting unit in forbidden unit.}}

\def\r@quired#1#2{%
\errhelp{Unit #2 must be inside unit #1, it isn´t.}
\errmessage{Starting unit outside of required unit.}}

\def\f@rbiddent@st#1\in#2#3{\def\\##1{\setr@g{t@st##1}\to\f@rbidden}%
\refreg{#1}%
\def\\##1{\e@\ifx\csname t@st##1\endcsname\f@rbidden\f@rbidden{##1}{#3}\fi}%
\reftokreg{#2}}

\def\r@quiredt@st#1\in#2#3{\def\\##1{\setr@g{t@st##1}\to\r@quired}%
\reftokreg{#2}%
\def\\##1{\e@\ifx\csname t@st##1\endcsname\r@quired\else\r@quired{##1}{#3}\fi}%
\refreg{#1}}

\def\forbidunit#1\in#2 {\leftappenditem #2\to f@rbiddens#1 }
\def\requireunit#1\inside#2 {\leftappenditem #2\to r@quireds#1 }
\def\forbidinself#1{\leftappenditem #1\to f@rbiddens#1 }
\def\norequireds#1{\e@\let\csname r@quireds#1\endcsname\relax}
\def\noforbiddens#1{\e@\let\csname f@rbiddens#1\endcsname\relax}

\def\namedef#1\in{\e@\langdef\csname#1name\endcsname\in}
\def\unitname#1{\csname#1 name\endcsname}
\def\presentunitname{\csname\presentunit name\endcsname}

\def\declareunit#1{%
\mkreg{r@quireds#1}\mkreg{f@rbiddens#1}%
\mktokreg{ev@ry#1}\mktokreg{@fter#1}}

\def\declareparent#1\of#2 {%
\e@\let\e@\s@tt\csname @ction#1\endcsname
\e@\let\csname @ction#2\endcsname\s@tt
\tokregconcat{ev@ry#2}={ev@ry#2}&{ev@ry#1} %
\tokregconcat{@fter#2}={@fter#1}&{@fter#2} }

\def\ignore#1{\everyunit{#1}{\afterunit{#1}{}}\unitaction{#1}{}%
\setr@g{ign@re#1}\to{\ign@re}}
\def\e@tgr#1{\egroup}
\def\ignorein#1{\everyunit{#1}{}\unitaction{#1}{\unittotoken{#1}\e@tgr}}

\long\def\everyunit#1#2{\csname ev@ry#1\endcsname{#2}}
\long\def\afterunit#1#2{\csname @fter#1\endcsname{#2}}

\long\def\addeveryunit#1#2{\tokleftappenditem{#2}\to{ev@ry#1} }
\long\def\addafterunit#1#2{\tokleftappenditem{#2}\to{@fter#1} }

\long\def\everyunitadd#1#2{\tokrightappenditem{#2}\to{ev@ry#1} }
\long\def\afterunitadd#1#2{\tokrightappenditem{#2}\to{@fter#1} }

\long\def\unitaction#1{\e@\def\csname @ction#1\endcsname}

\def\unittotoken#1#2{\t@={}%
\long\def\n@ste##1\end##2{\ifstreq{#1}{##2}%
\t@=\expandafter{\the\t@##1\end{##2}}#2\t@%
\else
\t@=\expandafter{\the\t@##1\end{##2}}%
\expandafter\n@ste\fi}\n@ste}

\let\@nd\end

\def\begin#1{\bgroup
\def\presentunit{#1}%
\r@quiredt@st{r@quireds#1}\in{presentpath}{#1}%
\f@rbiddent@st{f@rbiddens#1}\in{presentpath}{#1}%
\leftaddtoks\presentpath{\\{#1}}%
\setr@g{in@nit#1}\to{\in@nit}%
\def\sl@t{\reftokreg{ev@ry#1}\refreg{@ction#1}}%
{\ifunit{#1}%
 \aftergroup\sl@t
\else
\errhelp{The unit #1 has not been declared.}%
\errmessage{Undeclared unit}\fi}}

\def\end{\ifdef{presentunit}\e@\@nnd\else\e@\@nd\fi}
\def\@nnd#1{%
\let\sl@t\relax
\ifunit{#1}%
 \ifinunit{#1}\else
 \errhelp{The unit #1 has not been entered.}%
 \errmessage{Not in unit.}%
 \fi
 \ifpresentunit{#1}\else
  \edef\s@tt{The present unit is \presentunit,
  you are ending #1.}%
  \errhelp\e@{\s@tt}%
  \errmessage{Ending wrong unit}%
 \fi
\setr@g{in@nit#1}\to{\relax}
\def\sl@t{\reftokreg{@fter#1}}%
\else
\errhelp{The unit #1 has not been declared.}%
\errmessage{Undeclared unit}\fi
\sl@t\egroup}

\def\bye{\vfill \supereject \@nd}
\declareunit{document}
\everyunit{document}{\ignore{document}}
\afterunit{document}{\aftergroup\bye
\iffinal\else\medskip\parindent=0pt
\e@\if\the\x@bib\relax\relax
\else References:\par
\the\x@bib \medskip\fi\e@\if\the\x@note\relax\relax\else
Notes:\par\the\x@note\fi\fi\par}

\declareunit{start}
\everyunit{start}{\ifst@rted
\errhelp{You can only have one startunit in a document.}
\errmessage{Duplicate startunit.}
\else
\global\st@rtedtrue\fi}

\let\c@mma,
\catcode`\,=\active

\unitaction{start}#1\by{%
\count255=0
\@sfalse
\def\st@rtcent{\setbox0\hbox\bgroup\ignorespaces}
\def\@ndcent{\unskip\egroup\ifnum\count255>1%
\centerline{\unhbox2\tr@th\ \ifsvensk och\else and\fi\else\c@mma\fi}\fi
\advance\count255 by 1%
\tr@th\else\setbox2\box0\fi}
{\def\break{\hss\egroup\medbreak\hbox to\hsize\bgroup\hss}
	    \bf \centerline{#1}}%
                    \vskip6pt
                    \centerline{\ifsvensk av\else by\fi}
                    \vskip6pt
\vbox\bgroup\catcode`\,=\active
\sl
\def,{\@ndcent\par\st@rtcent}
\setbox0\hbox\bgroup\ignorespaces
}

\catcode`\,=12

\afterunit{start}{\global\@strue\@ndcent\par\centerline{\unhbox
0\unskip.}\egroup\vskip1in}

\declareunit{introduction}
\forbidinself{introduction}

\declareunit{chapter}
\everyunit{chapter}{\ch@pter}
\forbidinself{chapter}

\namedef appendices\in english={Appendix},\in
swedish={Appendix},\in french={Appendice}.
\declareunit{appendices}
\everyunit{appendices}{\vfill\eject
\noindent
\centerline{\bf\presentunitname}%
\nobreak\vskip1cm\nobreak\noindent}
\forbidinself{appendices}

\declareunit{appendix}
\unitaction{appendix}#1:#2\par{\vskip-\lastskip
\def\sektname{#1}%
\lop=0%
\bigbreak\noindent
\centerline{\bf#1: \ignorespaces#2}\nobreak\smallbreak\nobreak\noindent}
\requireunit appendix\inside{appendices}

\def\sectionhopp{\global\advance\sekt by 1}

\declareunit{sect}
\unitaction{sect}#1\par{\vskip-\lastskip
\global\advance\sekt by1%
\def\sektname{\number\sekt}%
\lop=0%
\bigbreak\noindent
\ifempty{#1}%
\sektname.\enspace
\else
{\bf\sektname. \ignorespaces#1}\par\nobreak\noindent\fi}

\declareunit{preliminaries}
\declareparent sect\of preliminaries
\addeveryunit{preliminaries}{\ifnum\sekt=0
\advance\sekt-1\else
\errhelp={Preliminaries are
only allowed as the first section of a chapter.}%
\errmessage{Preliminaries in the middle of chapter.}\fi}
\forbidinself{preliminaries}

\declareunit{section}
\declareparent sect\of section
\forbidinself{section}

\def\pr@{\par
 \ifdim\lastskip<\medskipamount \removelastskip\penalty55\medskip\fi}

\def\tr@mp#1.\@@#2{\def\@@{}\global\edef\t@mp{#1}}%
\declareunit{procl}
\unitaction{procl}#1 {%
\vskip-\lastskip
\global\ch@pqtrue
\global\advance\lop by 1%
\ifempty{#1}%
\edef\t@mp{\sektname.\number\lop.}%
\else
{\tr@mp#1\@@.\@@\@@}%
\e@\global\e@\edef
\csname l@p\t@mp\endcsname{\ifchap\the\chapno @\fi
\iffinal\sektname.\number\lop\else \t@mp\fi}%
\fi
\medbreak\noindent{\bf \presentunitname\
\iffinal
 \sektname.\number\lop.%
\else
 {\rm\ignorespaces\t@mp}
\fi}%
\enspace\sl}

\declareunit{stmt}
\declareparent procl\of stmt
\langdef\pr@@f\in english={PROOF},\in swedish={BEVIS},\in french={PREUVE}.
\addeveryunit{stmt}{\def\pro{\egroup\let\ifpr@\iffalse
\pr@{\bf \sevenrm\pr@@f}: }\bgroup\let\ifpr@\iftrue}
\afterunit{stmt}{\ifpr@
\errhelp={You seem to have forgotten to to finish the statement and
(at the same time) begin the proof with a \pro.}%
\errmessage{No proof in statement.}\egroup
\fi
\slut}

\declareunit{theorem}
\declareparent stmt\of theorem
\namedef theorem\in english={Theorem},\in swedish={Sats},\in
french={Th{\'e}or{\`e}me}.

\declareunit{proposition}
\declareparent stmt\of proposition
\namedef proposition\in english={Proposition},\in swedish={Proposition},%
\in french={Proposition}.

\declareunit{lemma}
\declareparent stmt\of lemma
\namedef lemma\in english={Lemma},\in swedish={Lemma},\in french={Lemme}.

\declareunit{definition-lemma}
\declareparent stmt\of definition-lemma
\namedef definition-lemma\in english={Definition-Lemma},%
\in swedish={Definition-Lemma},\in french={Definition-Lemme}.

\declareunit{corollary}
\declareparent stmt\of corollary
\namedef corollary\in english={Corollary},\in
swedish={F{\"o}ljdsats},\in french={Lemme}.

\declareunit{def}
\declareparent procl\of def
\afterunit{def}{\pr@}

\declareunit{definition}
\declareparent def\of definition
\namedef definition\in english={Definition},\in swedish={Definition},%
\in french={Definition}.

\declareunit{namedpar}
\everyunit{namedpar}{\vskip-\lastskip
\smallskip\noindent {\bf\presentunitname}\rm: \ignorespaces}
\afterunit{namedpar}{\par\smallskip}

\declareunit{remark}
\declareparent namedpar\of remark
\namedef remark\in english={Remark},\in
swedish={Anm{\"a}rkning},\in french={R{\'e}marque}.

\declareunit{example}
\declareparent namedpar\of example
\namedef example\in english={Example},\in
swedish={Exempel},\in french={Exemple}.

\declareunit{exercise}
\declareparent namedpar\of exercise
\namedef exercise\in english={Exercise \the\ovnnr},%
\in swedish={{\"O}vning \the\ovnnr},\in french={Exercise \the\ovnnr}.
\addafterunit{exercise}{\global\advance\ovnnr by 1}
\forbidinself{exercise}
\def\markexercise#1{\ifdef{ex@r@f#1}
\errhelp{#1 has already been used to mark an exercise.}%
\errmessage{Reuse of exercise mark.}%
\else
\e@\edef\csname ex@r@f#1\endcsname{\ovnnr}%
\fi}
\def\exerciseref#1{\ifdef{ex@r@f#1}\csname ex@r@f#1\endcsname
\else
\errhelp{#1 has not been used to mark any exercise.}%
\errmessage{Reference to unmarked exercise.}%
\fi}

\declareunit{solution}
\namedef solution\in english={Solution \the\ovnnr},%
\in swedish={L{\"o}sning \the\ovnnr},\in french={Solution \the\ovnnr}.
\requireunit solution\inside{exercise}
\forbidinself{solution}
\declareunit{hint}
\namedef hint\in english={Hint \the\ovnnr},%
\in swedish={Ledning \the\ovnnr},\in french={Suggestion \the\ovnnr}.
\requireunit hint\inside{exercise}
\forbidinself{hint}

\def\exerciselast{
\startout{exercise}\writeandput
\startout{hint}\writeonly
\startout{solution}\writeonly
\afterunitadd{document}{\ovnnr0
\unitaction{exercise}{}
\everyunitadd{exercise}{\ignorein{hint}}%
\everyunitadd{exercise}{\ignorein{solution}}%
\def\markexercise##1{}
\vfill\eject\centerline{\bf {\"O}vningar}
\vskip1cm
\csname input\endcsname\rjobname.exercise\
\declareparent namedpar\of hint
\norequireds{hint}
\vfill\eject\centerline{\bf Ledningar}
\vskip1cm
\csname input\endcsname\rjobname.hint\
\declareparent namedpar\of solution
\norequireds{solution}
\vfill\eject\centerline{\bf L{\"o}sningar}
\vskip1cm
\csname input\endcsname\rjobname.solution\
}}

{\catcode`\^^M=\active\catcode`\ =\active
\gdef\v@rb{\let^^M\ \let \ }
\gdef\v@rbdis{\let^^M\par\let \ }
\global\def\setesc@pe#1{\if#1 \catcode`@=0\else\catcode`#1=0\fi{}}
}

\declareunit{verb}
\everyunit{verb}{\bgroup\dospecials
\def\end{\egroup\end}\catcode`\^^M=\active\catcode`\ =\active\tt}
\unitaction{verb}{\setesc@pe}
\declareunit{verbatim}
\declareparent verb\of verbatim
\everyunitadd{verbatim}{\v@rb}
\forbidinself{verbatim}
\declareunit{verbatimdisplay}
\declareparent verb\of verbatimdisplay
\addeveryunit{verbatimdisplay}{\vskip-\lastskip\smallskip}
\everyunitadd{verbatimdisplay}{\v@rbdis\baselineskip10pt\parindent
0pt\leftskip2cm}
\afterunit{verbatimdisplay}{\vskip-\lastskip\smallskip\noindent}
\forbidinself{verbatimdisplay}

\declareunit{bibliography}
\everyunit{bibliography}{\bibl}

\def\final{\finaltrue\hfuzz=2pt\errorstopmode}
\def\pr@l@m{
}
\def\draft{\finalfalse\hfuzz=2pt\errorstopmode
\global\headline={%
\ifnum\pageno=1\hfill\pr@l@m
\else\tenrm\ifsvensk\idag\else\today\fi
 \hfill PRELIMINARY VERSION -- No preprint\hfill\rjobname\quad\fi}}

\def\t@stpt#1.#2.{#1\e@tif.}
\def\rjobname{\t@stpt\jobname..\relax}

\def\head#1\par{%
\global\headline={%
\iffinal
 {\ifnum\pageno>1\hfill\sevenrm
 \uppercase{#1}\hfill\else\hfill\fi}%
\else
 {\tenrm\ifsvensk\idag\else\today\fi
 \hfill\rjobname\quad}%
\fi}}

\def\ch@pter#1\par{
\ifch@pq
 \ifchap
 \else
  \errhelp={You must either use or not use chapters in a
document. In references
  chapternumbers will now be ignored.}%
  \errmessage{Chapter introduced in the middle of document.}%
 \fi
\else
 \global\chaptrue
\fi
\global\advance\chapno by1\sekt=0\lop=0
\vfill\eject
\noindent
\centerline{\bf\ifnum0=\chapno0\else
\uppercase\e@{\romannumeral\chapno}\fi: #1}%
\nobreak\vskip1cm\nobreak\noindent}

\def\fe@#1@{\ifnum#1=\chapno \else\ifnum0=#1 0\else
\uppercase\e@{\romannumeral #1}\fi:\fi}
\makesplit\spl@tco:

\def\isd@fin@d#1{%
\ifdef{l@p#1}%
 \@strue
\else
 \ifnoreferr
 \else
   \errhelp=\e@{#1 does not refer to anything.}%
   \errmessage{Reference to nonexistent result.}%
 \fi
  \@sfalse
\fi}

\def\(#1){%
(%
\spl@tco{#1}%
\edef\r@st{\csname l@p\firstsplit\endcsname\ \secondsplit\unskip)}%
\isd@fin@d\firstsplit
\ifchap
 \tr@th
   \e@\fe@\r@st
  \else
   #1)%
  \fi
\else
  \r@st
\fi
}

\def\ref#1{%
\isd@fin@d{#1}
\ifchap
 \tr@th
  \e@\e@\e@
  \fe@\csname l@p#1\endcsname
 \else(#1)%
 \fi
\else
 \csname l@p#1\endcsname
\fi
}

\def\slut{\nobreak\hfill \hbox{%
\vrule height 5pt
\kern-.4pt
 \vbox{
\hrule width 5pt depth0pt height.4pt
 \kern4.6pt \hrule  }
\kern-3.75pt
\vrule height 5pt}\kern1pt
\par}

\def\t@sta#1  #2{\t@{#1}\e@t}

\def\tag#1${%
\global\ch@pqtrue
\global\advance\lop by 1%
\t@sta#1 @ \unskip
\global%
 \e@\e@
 \e@\edef
 \e@\csname
 \e@
l\e@
@\e@
p\the\t@\endcsname
{
\ifchap
 \the\chapno @
\fi
\iffinal
 \sektname.\number\lop
\else
 \the\t@%
\fi
}
\tagchoice\hbox{(%
\iffinal
 \sektname.\number\lop
\else
 \the\t@%
\fi)%
}
$%
}

\def\tagchoice{\leqno}

\def\mkref#1{%
\global\advance\lop by 1%
\e@\global\e@\edef
\csname l@p#1\endcsname{\ifchap\the\chapno @\fi
\iffinal\sektname.\number\lop\else #1\fi}
\iffinal
 (\sektname.\number\lop)%
\else
 {\rm\ignorespaces(#1)}
\fi
}

\def\x@struf{\dp\strutbox}
\def\x@notnu{\strut\vadjust{\kern-\x@struf\vtop to\x@struf{%
\baselineskip\x@struf\vss\llap{\sevenrm\the\x@ant\quad}\null}}}
\x@sid={. s. }
\x@note={}
\long\def\note#1{%
\iffinal
\else
 \global\advance\x@ant1
 \iffinal
 \else
  \x@notnu
  \t@=
  \e@\e@
  \e@
  {\e@\the\e@\x@ant\e@\the\e@\x@sid
  \the\pageno: #1}%
  \global\x@note=\e@{\the\e@\x@note
  \e@\par\the\t@}%
 \fi
\fi}

\x@bib={}
\long\def\x@ny#1{\t@={#1}%
\global\x@bib=\e@{\the\e@\x@bib\e@\par\the\t@}}
\def\[#1]{[#1]\iffinal\else\x@ny{#1}\fi}
\namedef bibliography\in english={Bibliography},\in swedish={Bibliogafi},%
\in french={Bibliographie}.
\bibliographytrue
\def\bibl{
          \ifbibliography 
            \medbreak
            \centerline{\bf\presentunitname}
            \nobreak
            \medskip
            \nobreak
          \fi
          \ninerm
          \frenchspacing

\def\by{\egroup\setbox1=\hbox\bgroup\smc}
\def\paper{\egroup\setbox2=\hbox\bgroup}
\def\preprint{\egroup\prep@true\setbox2=\hbox\bgroup\sl}
\def\jour{\egroup\jour@true\setbox3=\hbox\bgroup\sl}
\def\vol{\egroup\setbox4=\hbox\bgroup\bf}
\def\yr{\egroup\setbox5=\hbox\bgroup}
\def\pages{\egroup\pages@true\setbox6=\hbox\bgroup}
\def\page{\egroup\page@true\setbox6=\hbox\bgroup p. }
\def\book{\egroup\book@true\setbox2=\hbox\bgroup\sl}
\def\lectnote{\egroup\lect@true\setbox7=\hbox\bgroup}
\def\publ{\egroup\setbox3=\hbox\bgroup}
\def\publadr{\egroup\setbox4=\hbox\bgroup}
\def\from{\egroup\setbox4=\hbox\bgroup}
\def\inbook{\egroup\inbook@true\setbox7=\hbox\bgroup}
\def\sln{\egroup\sln@true\setbox4=\hbox\bgroup}
\def\spec{\egroup\spec@true\setbox8=\hbox\bgroup}

\def\[##1:##2\par{\jour@false\book@false\inbook@false\pages@false
\sln@false\page@false\prep@false\spec@false\lect@false
\noindent\hangindent2.3\parindent\hangafter1\hbox
to2.3\parindent{[##1:\hfill}\bgroup##2\egroup\unskip
\unhbox1\unskip,
\ifjour@
\unhbox2\unskip, \unhbox3\unskip\ \unhbox4\unskip\
(\unhbox5\unskip), \unhbox6\unskip.\fi
\ifbook@
\unhbox2\unskip, \unhbox4\unskip, \unhbox3\unskip,
\unhbox5\unskip
\ifpages@\unskip, pp.
\unhbox6\unskip.\else\ifpage@\unskip
,\unhbox6\unskip.\else\unskip.\fi\fi\fi
\iflect@
\unhbox2\unskip, \unhbox7\unskip\ \unhbox4\unskip, \unhbox3\unskip,
\unhbox5\unskip
\ifpages@\unskip , pp.
\unhbox6\unskip.\else\ifpage@\unskip
,\unhbox6\unskip.\else\unskip.\fi\fi\fi
\ifinbook@
\unhbox2\unskip,
in \unhbox7\unskip, \unhbox4\unskip, \unhbox3\unskip,
\unhbox5\unskip, \ifpages@ pp. \fi  \unhbox6\unskip.\fi
\ifsln@
\unhbox2\unskip, SLN \unhbox4\unskip, Springer-Verlag
\ifpages@\unskip, pp.
\unhbox6\unskip.\else\ifpage@\unskip,
\unhbox6\unskip.\else\unskip.\fi\fi\fi
\ifprep@ \unhbox2\unskip, (preprint) \unhbox4\unskip, \unhbox5\unskip.\fi
\ifspec@\unhbox8\fi
\par}
}

\def\sc@le#1#2{\divide#1 by 1000\multiply#1#2}

\def\p@llbox#1#2#3{\hbox{%
\getsc@le{#1}%
#2%
\ifx\scalep@rt\relax\def\scalep@rt{1000}\else
\sc@le\h@ight\scalep@rt
\sc@le\w@dth\scalep@rt\fi
\vbox to \h@ight{\hrule width \w@dth height 0pt depth 0pt
#3{#1}%
\vfill\special{\t@mp}}}}

\def\eatsc@led scaled{}
\def\insc@le#1scaled#2scaled{\global\def\namep@rt{#1 }
\ifempty{#2}\global\let\scalep@rt\relax\else
\global\def\scalep@rt{#2}\e@\eatsc@led\fi}

\def\getsc@le#1{\insc@le#1scaledscaled}

\newread\@psffile
\dimendef\h@ight=0
\dimendef\w@dth=1

\def\pictg@tdim#1#2{\h@ight=#2\unskip\w@dth=#1\unskip}

\def\pictt@mpl#1{\def\t@mp{picture #1}}
\def\pictbox#1 by #2 (#3){\p@llbox{#3}{\pictg@tdim{#1}{#2}\unskip}\pictt@mpl}

\def\picture #1 by #2 (#3){\vskip 6pt\centerline{%
\pictbox #1 by #2 (#3)}\vskip 6pt\noindent}

{
\catcode`\%=11
\catcode`\|=14

\global\def\illg@tdim{|
\openin\@psffile=\namep@rt \unskip
{|
\catcode`\%=11|
\catcode`\_=11|
\catcode`\&=11|
\def\h@ndle##1 ##2 ##3 ##4:{{\globaldefs1\def\l@ft{##1 bp}|
\def\b@t{##2 bp}\def\r@ght{##3 bp}\def\t@p{##4 bp}}}|
\def\t@est{
\def\e@tspace##1{\if##1 \else##1\fi}|
\def\BBt@st##1:##2
\edef\s@tt{\e@tspace##2}|
\e@\h@ndle\s@tt\unskip\let\n@xt\relax\fi}|
\read\@psffile to \h@p
\def\n@xt{|
\read\@psffile to \h@p
\e@\BBt@st\h@p :
\n@xt}|
\n@xt
}|
\closein\@psffile
\e@\h@ight\t@p
\e@\advance\e@\h@ight\e@-\b@t
\e@\w@dth\r@ght
\e@\advance\e@\w@dth\e@-\l@ft}
}

\def\illt@mpl#1{\def\t@mp{illustration #1}}
\def\illbox(#1){\p@llbox{#1}\illg@tdim\illt@mpl}

\def\illustration(#1){\vskip 6pt plus 2pt minus 2pt\centerline{\illbox(#1)}%
\vskip 6pt plus 2pt minus 2pt\noindent}

\let\psprol\relax
\let\psend\relax
\def\pst@mpl#1{\e@\w@dth\l@ft \e@\h@ight\b@t
\sc@le\w@dth\scalep@rt \sc@le\h@ight\scalep@rt
\divide\w@dth65536 \divide\h@ight65536
{\count255=\scalep@rt \divide\count255 10
\count254\w@dth \count253\h@ight
\xdef\t@mp{psfile=\namep@rt hoffset=-\number\count254 \space\space
voffset=-\number\count253\space\space vscale=\number\count255\space\space
hscale=\number\count255}%
}}

\def\psbox(#1){\psprol\p@llbox{#1}\illg@tdim\pst@mpl\psend}

\def\psillustration(#1){\vskip 6pt plus 2pt minus 2pt\centerline{\psbox(#1)}%
\vskip 6pt plus 2pt minus 2pt\noindent}

\ifmac\else
\let\illbox\psbox
\let\illustration\psillustration
\fi

\long\def\condition#1#2\endcondition{\medskip\setbox14=\hbox{#1
}\setbox13=\vbox{\advance\hsize-1.5cm\advance\hsize-\wd14\noindent #2}
\dimen1=\ht13 
\hbox to\hsize{\hfill\vbox
to\dimen1{\vfill\box14\vfill}\hfill \box13\hskip.75cm}\medskip\noindent}

\langdef\yrs\in english={Yours\iffriendly\else sincerely\fi},%
\in swedish={\iffriendly H\else V{\"a}nliga h\fi {\"a}lsningar},%
\in french={Bien á \iffriendly toi\else vous\fi}.

\def\Torsten{\vskip1cm\leftline{\hskip10cm \yrs}\nobreak \vskip.75cm\nobreak
\leftline{\hskip10cm Torsten \iffriendly\else Ekedahl\fi}}

\def\letter{\footline={\hss\ifnum\pageno=1\else\tenrm\folio\fi\hss}%
\headline={\ifnum\pageno>1T.~Ekedahl\hss\today\else\hfil\fi}%
\rightline{Stockholm \today\qquad}\vskip3cm}

\langdef\dear\in english={Dear},\in swedish={},\in french={Cher}.
\def\beginletter#1\par{\letter \dear #1\unskip,\par\qquad}

\def\endletter#1\par{\if\relax#1\relax
\else\friendlytrue\fi\nobreak\Torsten\vfill\eject\@nd}

\catcode`\@=12

\ifdump
\let\dmp\dump
\else
\let\dmp\relax
\fi
\dmp

\final
\mdef\barfun{\overline{\hbox to 5pt{\vrule width 0pt height 5pt\hss}}}

\begin{document}

\begin{start}
Boundary behaviour of Hurwitz schemes
\by
Torsten Ekedahl
\end{start}
\begin{introduction}
The purpose of this paper is to get a convenient way of computing the
topological types of degenerations of covers of the projective line. This leads
immediately to the theory of Hurwitz schemes.
Hurwitz schemes are normally thought of as the moduli spaces
of maps of curves to
$\P^1$ with prescribed degree and local behaviour. We will rather consider the
Galois hull of the mapping and thus consider curves together
with an action of a
finite group with only trivial invariants on 1-cohomology. If the group then is
provided with a transitive permutation representation, a curve with a map to
$\P^1$ is obtained and the local behaviour is easily described in terms of the
(conjugacy classes of) stabiliser subgroups of points. Now, if the curve and
the group action varies in a family the family of curves has, after a
generically finite base extension, a stable extension to a proper base and the
group action automatically extends. If we look at boundary points we will
then get an action of our group $G$ on a stable curve. What
this note actually does
is to provide a description of such actions in terms of homomorphisms from a
kind of fundamental group of the quotient curve into $G$, a
description which is
completely analogous to the well known case of a smooth curve. A direct
compactification of the Hurwitz scheme was provided by Harris and Mumford in
\[Ha-Mu] and it is clear that one could extract a description of the possible
topological types from their construction. The present paper differs from such
an approach mainly in emphasis. Here we will take the point of view of actions
of finite groups on stable curves, which leads, I believe, to an easier way of
understanding the topological types of the degenerations. If one tries to make
the approach chosen in this paper into a definition of a
moduli problem one will
also get a problem slightly different from the one of Harris and Mumford,
showing that the two possible approaches are not formally equivalent. The
moduli point of view will however not be important to us. We will in this paper
take the algebraic point of view. Over the complex numbers an analytic approach
would, of course, be just as possible. We will from time to time indicate the
changees necessary in an analytic setting. The same goes for
the change from the
view point of stacks, which is that taken here, to orbifolds.
\end{introduction}
\begin{preliminaries}Preliminaries

We will need to compare the cohomology of a curve and a stable reduction of
it. Normally this would be done using vanishing cycles.
However, here we will be
only interested in rather crude properties of cohomology (essentially its
character under a finite group action) so we will give a,
for us, more convenient
description.
\begin{remark}
The end result is of course well known irrespective of if one uses vanishing
cycles or not, only that we need to be careful about identifications as we will
have a group acting on the curve and thus on cohomology.
\end{remark}
Recall that for a curve $C$ (which to us will be reduced and proper over an
algebraically closed field) we have both the sheaf of K{\"a}hler 1-forms
$\Omega^1_C$ and the dualising sheaf $\omega_C$. We further have a trace map
$\Omega^1_C \rightarrow \omega_C$\note{Referens}. Composing this trace map with
the total differential we get a complex of sheaves $\Cal O_C\rightarrow
\omega_C$, which we will call $\omega^._C$. As $C$ is
Cohen-Macaulay it is clear
that the cohomology of this complex is self-dual just as in the smooth
case. Furthermore, the spectral sequence abutting to the
hypercohomology of this
complex degenerates for obvious reasons showing that
$H^0(C,\omega^._C)=H^0(C,\Cal O_C)$, $H^1(C,\omega^._C)$ has $H^0(C,\omega_C)$
as a subspace with quotient $H^1(\Cal O_C)$ and
$H^2(C,\omega^._C)=H^1(C,\omega_C)$. We will, for obvious reasons, call this
hypercohomology the {\deffont de-Rham cohomology} of $C$.

In particular if $C$ varies in a flat family of curves, then we get a relative
complex $\omega^._{C/S}$ whose higher direct images to $S$ will be flat sheaves
commuting with base change. The descriptions of the hypercohomology given will
also be true in a family.
\begin{remark}
i) Note that in general these higher direct images will not have integrable
connections, this seems to be the major difference from the smooth case. We
will only be using the fact that the character of a finite group is constant in
a (connected) family and we do not need an integrable connection for that.

ii) As we are really only interested in the character on
cohomology for a smooth
curve described in terms of the action of the group on a degeneration we could
also have argued that the virtual character on the smooth curve is the virtual
character on the cohomology of the structure sheaf plus its complex conjugate
and then only used that the cohomology of the structure sheaf behaves well in
families. The definition given here of de Rham-cohomology could however be of
independent interest. (When all singularties are nodes it is the log-de Rham
complex associated to a log-structure of the curve.)
\end{remark}
Let us now concentrate on the case which will interest us, namely when $C$ has
only ordinary double points. Let us agree to use the term
\definition {branch at
a node} to denote one of the points above the node in the normalisation of the
curve (the standard terminology being to let it denote one of local irreducible
components at the node, there is a natural and canonical
bijection between these
and our branches). We will then compare the cohomology of
$\omega^.$ for $C$ and
its normalisation \pil\pi{C'}C. To begin with we have a map $\Cal
O_C\rightarrow\Cal O_{C'}$ whose cokernel is the sum over the double points of
$C$ of copies of the base field. Note, however, that that copy is not
canonically isomorphic to the base field but rather isomorphic to the free
vector space on the two branches divided by the sum of the basis elements. Thus
an automorphism that switches the two branches acts on the copy by
multiplication by $-1$. Dualising this we get a map
$\omega_{C'}\rightarrow\omega_C$. The cokernel may again be identified with the
2-dimensional vector space based on the branches divided by the sum of the base
vectors. Taking long exact sequences of cohomology gives a subspace of
$H^1(C,\Cal O_C)$ isomorphic to the cokernel of the quotient map $\Cal
O_{C'}\rightarrow\Cal O_{C'}/\Cal O_C$ whose quotient is isomorphic to
$H^1(C',\Cal O_{C'})$. The description of $H^0(C,\omega_C)$ is dual.
\begin{remark}
In the case of the cokernel of $\omega_{C'}\rightarrow\omega_C$ the map to the
1-dimensional vector space may be seen clearly by
identifying the sheaf pullback
of $\omega_C$ to $C'$ with 1-forms on $C'$ with at most simple poles at the
branches and for which the sum of the residues at the branches of a fibre is
zero. This also makes it clear that an automorphism switching the two branches
acts by multiplication by $-1$.
\end{remark}
We are now almost at the point where we can summarise this d{\'e}vissage of the
cohomology. However, as has already been mentioned, we will be interested only
in the action of a finite group and in fact, as far as cohomology is concerned
we will only care about the virtual character for the group. It is clear that
the dual graph and its cohomology should appear in that d{\'e}vissage. As may
already be apparent from the discussion above some care has to be taken in how
to define the dual graph. We will need to adopt the definition of \[Se:2.1.1]
of a graph.
\begin{definition}
i) A {\deffont graph} consists of a {\deffont vertex set} $V$, a set $V$ of
{\deffont oriented edges}, a map
$$
\funk V{E\times E}v{(o(v),t(v)),}
$$
 and an involution \pil\barfun VV such that $o(\overline v)=t(v)$ and
$\overline v\ne v$ for all vertices $v$. The vertices $o(v)$ and $t(v)$ will be
called the {\deffont source} and the {\deffont sink} of $v$ respectively. The
equivalence classes under \barfun{} of $V$ will be called the {\deffont edges}
of the graph. Automorphisms of a graph are permutations of
$E$ and $V$ commuting
with $o$, $t$ and \barfun. If $v$ is an oriented edge $\overline v$ will be
called its {\deffont opposite} edge.

ii) A \definition {generalised graph} consists of the same data and conditions
as a graph with the exception that \barfun{} is not assumed
to be fixpoint free.
\end{definition}
\begin{remark}
Our graphs will mainly be dual graphs of stable curves which
will always fulfill
the conditions for the definition of a graph. However, when a finite group is
acting on a stable curve we will be interested in the
quotient of its dual graph
by the group and such a quotient will, in general, only be a generalised
graph. The reason why we make the distinction is that the discussion of the
1-chains would not be completely true for a generalised graph (an edge equal to
its own opposite would give an element of order 2 in the group of 1-chains).
\end{remark}
The chain complex of a graph is defined as follows. We let the 0-chains be the
ordinary 0-chains but let the 1-chains be the free abelian group on oriented
1-simplices divided by the subgroup generatad by the elements which are the sum
of an edge and its opposite. Note that this 1-chain group is isomorphic to the
free group on the edges of the graph but not canonically so (or maybe
canonically but not naturally isomorphic). In particular an automorphism that
fixes an edge but permutes the oriented edges of that edge will act by
multiplication by $-1$ on the corresponding element of the
1-chains. Representation-theoretically, if a group acts on
the graph we see that
the 0-chains are isomorphic to the permutation representation on vertices,
whereas the representation on 1-chains is the sum over the orbits of the group
on edges of the representations induced from signum characters. By the latter I
mean the following. Pick a representing edge of an orbit and consider its
stabiliser. Each element of the stabiliser will either fix the two oriented
edges or permute them. If we assign $1$ resp.~$-1$ to the element we get a
character of the stabiliser which we will call the signum
character. We can then
induce this character to a representation of the whole group. We will call an
orbit where the signum character is trivial an {\deffont orientable} orbit and
in the opposite case we will speak about an {\deffont unorientable} orbit, in
the orientable case the representation is of course just the permutation
representation.  Note now that if the group acts on the curve then that an
automorphism fixes an edge means exactly that it fixes the corresponding double
point and that it fixes the edge's two vertices means that it fixes the two
branches. Thus we see that the space of global sections of $\Cal O_{C'}$ is
(canonically) isomorphic with the 0-chains (more naturally 0-cochains) with
scalars extended to the base field and that the global sections of $\pi_*\Cal
O_{C'}/\Cal O_C$ is (canonically) isomorphic with the 1-chains again with
scalars extended to the base field. A moment's thought will also convince the
reader (though we will not need it) that the map induced on
global sections from
the quotient map equals the coboundary map on cochains.
Hence by the d{\'e}vissage
above we see that $H^1(C,\Cal O_C)$ contains a subspace isomorphic to the first
cohomology space of the dual graph (with coefficients in the base field) with
quotient space equal to $H^1(C',\Cal O_{C'})$. There is also a dual description
for $H^0(C,\omega_C)$. We summarise this in the following, where by the
\definition{virtual character} of some cohomology we, naturally, mean the
alternating sum of the character of the individual cohomology groups.
\begin{proposition}
Let $C$ be a curve with only ordinary double points and $C'$ its
normalisation. Assume $C$ is acted upon by a finite group $G$.

i) The virtual character of the cohomology on the de-Rham cohomology of $C$ is
equal to the sum of that of the de-Rham cohomology of $C'$ and 2 times
the virtual character of the cohomology of the dual graph.

ii) The virtual character of the cohomology of the dual graph is equal to the
difference of the permutation character of the vertices minus the sum over the
orbits of $G$ on the edges of the representations induced from the signum
representation.
\pro The two statements follow directly from the previous discussion.
\end{proposition}
When we later want to discuss the description of a stable curve with an action
by a group $G$ we will need to describe this fundamental group in terms of
graphs of groups (cf.~\[Se:4.4.8]). However, it turns out that the
notion we need is a slight extension of this notion to {\it
generalised} graphs.
\begin{definition}
i) A (generalised) \definition {graph of groups} consists of
a generalised graph
and the association of a group $G_v$ to each vertex of the graph, a group $H_y$
to each edge of the graph, a group homomorphism \pil{a\mapsto
a^y}{H_y}{G_{t(y)}} for each edge $y$ and an isomorphism
\pil\barfun{G_y}{G_{\overline y}} whose square is the identity.

ii) The \definition {fundamental group} of a generalised graph of groups is the
same as that of \[Se:p.~42] with the modification that (using the
notation of (loc.~cit.)) $ya^y\bar y=\bar a^{\bar y}$.
\end{definition}
\begin{remark}
When the generalised graph is a graph one can use the isomorphisms $G_y\to
G_{\overline y}$ to identify $G_y$ and $G_{\overline y}$ and then we
arrive at the notion defined by Serre (except that Serre assumes that $a\mapsto
a^y$ is injective, on the one hand that condition will not be necessary for our
purposes, on the other hand it will always be fulfilled in our applications).
\end{remark}
The more precise way in which graphs of groups come up in our calculations is
through the use of the van Kampen theorem. As we will have occasion to apply
this theorem to the case where the intersection is not connected, we will
quickly recall the results in that case as it may be less
known than the case of
connected intersection. We are also forced to formulate the
result for algebraic
stacks (orbifolds) as we are going to apply it in that generality.
\begin{definition-lemma}vK
Let $U_1$ and $U_2$ be an open covering of a connected algebraic stack
$X$. Assume that an {\'e}tale cover of $U_1$, one of $U_2$ together with an
isomorphism of their respective restrictions to the fibre product of $U_1$ and
$U_2$ over $X$ comes from a unique cover $X$. Then the
fundamental group of $X$ (wrt some base point) is isomorphic to the fundamental
group of the following graph of groups:

Take one vertex for each connected component of $U_1$ and
one for each connected
component of $U_2$. Attach to these vertices the fundamental group of the
corresponding component. For each connected component $V$ of
$U_1\bigcap U_2$ and
components $U'_i$ of $U_i$, $i=1,2$ which both contain $V$
one defines a pair of
oriented arrows from $U'_1$ to $U'_2$ and vice versa. Attach to these edges the
fundamental group of $V$ and the morphism from that fundamental group to the
fundamental group of $U'_i$ induced by inclusion.
\begin{remark}
The condition would seem to be always verified. However,
note that, by necessity
an open covering means an open {\'e}tale cover (as we are dealing with
stacks). Hence one would in general also need descent conditions on the
``intersection'' (i.e.~fibre product) of the $U_i$ with themselves.
\end{remark}
\pro The proof becomes more understandable if one passes to fundamental
groupoids instead of fundamental groups (which also explains why base points do
not play any role in the formulation of the result). To do that we define the
notion of map from a graph of groups to a groupoid (in analogy with the
definition of \[Se:1.1.1]) as a map taking vertices to objects, a
homomorphism from the group at a vertex to the automorphism group of the
corresponding object and a map from edges to morphisms the source object to the
sink object, fulfilling the obvious relations. Two maps which are related by an
equivalence of groupoids are then appropriately identified
so as to give rise to
the notion of 2-limit of a graph of groups. It is then an easy exercise to show
that the fundamental group of a graph of groups is a 2-limit in the category of
groupoids. The point of considering groupoids instead of groups is that the
proper formulation of the van Kampen theorem when the intersection is
non-connected is that the fundamental groupoid of the union is the 2-limit of
the fundamental groupoid of the parts with gluing of the intersection (this
follows for instance directly from the description of the fundamental groupoid
as constructed from the category of {\'e}tale coverings).
Gathering together these
strands gives the result.
\end{definition-lemma}
\end{preliminaries}
\begin{section}Group quotients of stable curves.

As our first step we will study the quotient of a stable
curve by a finite group
of automorphisms and then discuss what topological data are needed to recover
the original curve from its quotient. We begin by recalling the following
observation from \[Ka-Ma] which shows that we do not have to
worry about whether
we want to take a quotient of a family or a member of that family.
\begin{lemma}
Let $X\to S$ be a map and $G$ a finite group acting on $X$ over $S$. If the
order of $G$ is invertible in $\ko S$ taking the quotient by $G$ commutes with
base change.
\pro See \[Ka-Ma:A7.1.3.4].
\end{lemma}
\begin{remark}
We will only be interested in the case when the order of the group is indeed
invertible so that from now on {\it the order of a group
acting on a curve will,
unless otherwise mentioned, be supposed to be invertible in the base}.
\end{remark}
This means that in order to study the properties of a quotient of a stable
family of curves we should often be able to reduce to the base being a field.
Note that if we have a curve $C$ with a finite group $G$
acting on it and $x$ is
an ordinary double point on $C$ then the stabiliser has the following
structure. First there is a subgroup, $G'_x$, of index at
most 2 of the elements
which preserve the two branches. That subgroup acts on the
tangent spaces of the
two points of the normalisation lying over $x$, giving rise
to two characters of
that group. We will call these two characters the {\deffont characters of the
branches}.
\begin{definition}
Let $C$ be a curve, $G$ a finite group acting on it and $x$ an ordinary double
point of $C$. The action of $G$ is said to be {\deffont admissible} at $x$ if
the two character of the branches are inverses to each other and the square of
any element of the stabiliser $G_x$ of $x$ exchanging the two branches at $x$
acts trivially in a neighbourhood of $x$.
\end{definition}
\begin{remark}
We will see later that the action is admissible at a point iff that point
is equivariantly smoothable.
\end{remark}
It will be convenient to have the following description of
admissibility which is
almost purely group theoretical.
\begin{definition-lemma}
Let $C$ be a curve, $G$ a finite group acting on it and $x$ an ordinary double
point of $C$. Let $H_x$ be the stabiliser of $x$ modulo those elements which
act trivially on a neighbourhood of $x$. Then the action of
$G$ is admissible at
$x$ iff either

i) $H_x$ is dihedral of order greater than 4

\noindent or

$H_x$ is dihedral of order 4 or 2 with some
element exchanging branches and the characters of the branches being equal.

\noindent or

ii) $H_x$ is cyclic and the two characters of the branches are each other's
inverses.

\noindent In these two cases we will say that $x$ is a
\definition {dihedral} resp.~\definition {cyclic} node (for the given action of
$G$).
\begin{remark}
The formulation of the condition in i) is of course
complicated by the fact that
the ``dihedral'' groups of order 2 and 4 differ from the other dihedral groups
in that they are not really identifiable as such.
\end{remark}
\pro If there are no elements of $H_x$ exchanging the branches then the two
characters give a faithful representation of the whole group. If they are
inverses to each other then that means that the group is cyclic. This leads to
case ii) where the converse is obvious. If there is an \gsi{} element of $H_x$
exchanging the branches then such an element has order 2 (if the point is
admissible) and the subgroup $H'_x$ of index 2 of $H_x$ is cyclic as in the
first part. Furthermore, \gsi{} exchanges the branches and thus exchanges their
characters. As they are faithful and inverses to each other conjugation by
{\gsi} has to invert any element of $H'_x$ and thus $H_x$ is
dihedral. Conversely, where we may assume that the order of $H_x$ is greater
than 4, if $H_x$ is dihedral then there must be an element exchanging the
branches because the subgroup of $H_x$ fixing the branches is abelian. For the
same reason this element must be outside the cyclic subgroup of index 2. As
conjugation by any such element acts by multiplication by $-1$ on the cyclic
subgroup, the two characters of the branches are inverses to each other.
\end{definition-lemma}
The first structure result on quotients that we will obtain is a description of
the singularities of it; as we will see all singularities will again be
nodes. Note that as we have already concluded that quotients commute with base
change the conclusion is that if we have a family of curves all of whose
singularities are nodes then the same is true of a quotient by a finite group
acting admissibly.
\begin{proposition}
Let $C$ be a curve all of whose singularities are nodes and
$G$ a finite group of
automorphisms of $C$ such that the action is admissible at
all nodes. Then $C/G$
has only nodes as singularities and the singularities are the points below the
cyclic nodes of $C$.
\pro The quotient of a smooth point is normal and thus again smooth. At a node
we may localise and complete and then, as the order of $G$ is invertible in the
base field, we may choose a $G_x$-invariant complement to the square of the
maximal ideal in the maximal ideal; $m=V\Dsum m^2$. When $x$ is cyclic $V$ is
the direct sum of the two characters of the branches and
when $V$ is dihedral it
is irreducible. In any case we may choose an isomorphism of the completed local
ring with $\pow{x,y}/(xy)$ in a way such that up to multiplication by scalars
$x$ and $y$ are either fixed or permuted by elements of $G_x$. In the cyclic
case they are always fixed and if the order of the
characters of the branches is
$n$ then the fixed ring is $\pow{x^n,y^n}/(x^ny^n)$ which shows that the
quotient has a nodal singularity. in the dihedral case we may first take
invariants under the  subgroup of index 2. As we have seen that fixed ring
has a nodal singularity so that we are reduced to the case when the $G_x$ has
order two and the involution permutes $x$ and $y$. Then the ring of invariants
is simply $\pow{x+y}$ which is a regular ring.
\end{proposition}
If we now want to emulate the description of actions of groups on smooth curves
in terms of the quotient and topological data on it to reconstruct the ramified
cover, we are faced with the problem that quotients of dihedral points are
smooth yet we evidently need to treat them differently from those lying below
smooth points. What is happening here is that the quotient should really
be thought of as an algebraic stack with the points below dihedral points being
non-scheme points. Furthermore, at a singular point of the base, a quotient map
will in general not be {\'e}tale nor will any map which is
{\'e}tale in a punctured
neighbourhood be accepted. Again the necessary restrictions on the local
behaviour can be formulated in terms of gluing in a non-scheme point at a
singular point and as we have already decided to take the plunge into algebraic
stacks it seems reasonable to adopt that approach for those points as
well. Finally, for reasons of symmetry it seems reasonable to treat ordinary
ramification points in the same manner though that is clearly not traditional,
nor is it necessary.
\begin{remark}
i) As an apology to readers who feel ill at ease with algebraic stacks I would
offer that I started with the ambition to simplify their lives by avoiding the
use of stacks but eventually found that approach too cumbersome to feel able to
uphold that ambition.

ii) Depending on one's choice of education one may, of course, be used to talk
of orbifolds rather than algebraic stacks. For readers with this tendency I
would like to point out a few differences.
\item{a)} Usually orbifolds are assumed to be smooth, that
is, modeled on actions
by groups on smooth manifolds. We will have occasion to consider those which
are not (and thus maybe should be called ``orbispaces'' instead).
\item{b)} Usually, orbifolds are modeled on quotients by finite groups, here we
will need to use infinite groups.
\item{c)}Here we are not only considering things topologically but also
algebraically-analytically, in fact we will really only do things algebraically
even though the changes to make it work in the complex analytic category are
trivial.
\end{remark}
To prepare for the definition of the appropriate algebraic
stacks we will define
how to glue in the desired special points. As we are only considering group
actions where the order of the group is invertible in the base field we will in
the same vein only consider tame coverings. Let $C$ be a curve and $c\in C$ a
smooth point on $C$. Let $R$ be the (strict) Henselisation of the local ring of
$C$ at $c$. Let $\tilde R$ be the extension of $R$ obtained
by adjoining all the
roots, of order prime to the characteristic of the base field, of a uniformiser
of $R$ (this of course is the maximal tamely ramified extension of $R$). The
disjoint union, $D$, of $\Sp R$ and $C\setminus {c}$ maps to
$C$. If we take the
fibre product over $C$ of this disjoint union then outside
of $c$ it is an {\'e}tale
equivalence relation. We will modify that fibre product at $c$ so that it
together with $D$ forms an {\'e}tale groupoid. In fact, $D$
has a component $\tilde
P$ which is the fibre product of $\Sp \tilde R$ with itself over
$\ko{C,c}$. This in turn is a scheme over the fibre product of $\Sp R$ with
itself over $\ko{C,c}$. This latter scheme is the disjoint union of schemes
some of which maps isomorphically to $\Sp R$ through any of
the projections (and
which are parametrised by the Galois group of the residue field of $\ko{C,c}$)
and spectra of fields. As these latter maps into $C\setminus {c}$ we may ignore
them wanting only to modify the fibre product at $c$. For
each of the components
of the former type, the inverse image of $\tilde P$ over it
is isomorphic to the
fibre product $Q$ with itself of $\Sp \tilde R$ over $\Sp R$ and so we may
consider instead that case and then simply transfer the
modifications to $\tilde
P$. Now, we have an action of the Galois group \Cal G of the fraction field of
$\tilde R$ over $R$ (which is isomorphic to $\hatt \Z$ with
the component $\Z_p$
removed when the base field has positive characteristic) and
hence we have a map
$\Cal G\times \Sp{\tilde R} \to Q$. This is an isomorphism over the fraction
field of $R$ and our modification consists exactly of replacing $Q$ by $\Cal
G\times \Sp{\tilde R}$. That we in this fashion do indeed
get an {\'e}tale groupoid
is easily checked. This defines an algebraic stack $C_c^{cycl}$. It has an
open subset isomorphic to $C\setminus {c}$ and the complement consists of a
closed point. We will call this process {\deffont gluing in a cyclic point at
$c$}. It can evidently be repeated so that we can glue in cyclic points at a
finite set of points.
\begin{remark}
Intuitively, we have taken a small neighbourhood of $c$, taken the maximal
tamely ramified cover of the punctured neighbourhood, added a point to that
cover lying over $c$ and then glued in the stack quotient by
the Galois group of
the cover, where the ordinary quotient is just the original neighbourhood. That
intuitive picture was somewhat hidden both by not making the
assumption that the
residue field of $c$ should be algebraically closed and by the fact that an
{\'e}tale neighbourhood is not a subset. If one works in an analytic setting
instead then the truth is indeed very close to this intuition. Indeed, adding a
point to the universal cover of a punctured neighbourhood of $c$ is exactly the
procedure of adding a cusp at $i\infty$ to the upper half
plane. Thus taking the
stack quotient of the upper half plane with cusps added at rational points of
$\P^1(\R)$ by $PSL(2,\Z)$ is the same thing as compactifying
the moduli stack of
elliptic curves and then gluing in a cyclic point at $\infty$.
\end{remark}
We have also set up things so that a covering $D\to C$ which is tamely ramified
over $c$ extends uniquely to a cover of $C_c^{cycl}$ which
is {\it {\'e}tale} over
$c$ and, conversely, a cover of $C_c^{cycl}$ {\'e}tale over
$c$ gives a covering
$D\to C$ which is tamely ramified over $c$.

We will now slightly modify this construction so as to make it apply in the
other two situations mentioned above. First we consider
again a curve $C$ with a
smooth point $c$. We may again also consider $R$, the strict henselisation of
the local ring at $c$, and $\tilde R$, the maximal tamely ramified extension of
$R$. We now consider the direct product of two copies of $\tilde R$ as a ring
and the subring of that product consisting of the pairs $(r,s)$ for which $r$
and $s$ are congruent modulo the maximal ideal of $\tilde R$. This new ring
$\hatt R$ has an action of the Galois group of (the fraction field of) $\tilde
R$ over (the fraction field of) $R$, which acts the same way on both factors of
the direct product. We also have the automorphism of $\hatt R$ which exchanges
the two factors. Together with the automorphisms coming from
the Galois group it
generates a dihedral group. We now continue the construction exactly as in the
previous case using $\hatt R$ instead of $\tilde R$. We will call this process
{\deffont gluing in a dihedral point at $c$}. It can
evidently be repeated so that
we can glue in dihedral points at a finite set of points (as well as some
dihedral and some cyclic points).

Finally, we will describe how to glue in a cyclic point at an ordinary
node.  Let us first recall some facts
about algebraic fundamental groups and ``paths''. Thus if $R$ is a discrete
valuation ring the tame fundamental group of its fraction field is canonically
identified with $\hat\Z(-1) [1/p]$, the dual of the Tate module of roots of
unity in an algebraic closure of the fraction field; $p$
being the characteristic
exponent of the fraction field. In particular this group only depends on the
choice of an algebraic closure of the residue field of $R$. This means that if
we have an ordinary double point and if we consider the (strict) Henselisation
of of its local then ring, then the tame fundamental group of its two fraction
fields are canonically identified.
\begin{remark}
We haven't mentioned  complex analytic analogues of our
constructions but they are
clearly possible. In the analytic case a choice of
orientation of {\C} identifies the
local fundamental group with {\Z} and thus the fundamental groups of the two
components of a small punctured neighbourhood of an ordinary double point are
also canonically identified and the gluing in of a cyclic point at an ordinary
node can be done.
\end{remark}
 We start now with a point $c$ which is supposed to be an ordinary
double point. Again we consider $R$, the henselisation of the local ring. Its
total ring of fraction is the product of two fields. For each of those we take
the maximal tamely ramified extension. The normalisation of $R$ in the product
of these two fields is itself the product of two copies of the same local
ring. We then let $\tilde R$ be the subring of that ring of pairs $(r,s)$ for
which $r$ is congruent to $s$ modulo the maximal ideal. Again the Galois group
of the maximal tamely ramified extension acts on $\tilde R$, where we let the
action by given by the natural action on one of the copies and the inverse of
the natural action on the other copy.  We then continue as before.
\begin{remark}
As was mentioned above these constructions are introduced to describe certain
coverings. If we repeat this motivation together with some more concrete facts
we get the following.

i) For a smooth cyclic point any covering of the curve extends uniquely to a
covering of the stack which is {\it {\'e}tale} at the cyclic point.

ii) For a cyclic ordinary double point any covering of the curve outside of the
point extends uniquely to a covering of the stack which is
{\it {\'e}tale} at the
cyclic point iff the ramification indices at the two branches are the
same. Hence for instance $s^2=x,t^3=y$ does not so extend over $xy=0$.

iii) For a dihedral point the situation is different. Just as in case ii) a
covering may not extend to an {\'e}tale covering of the
stack. The main difference is
that the extension, if it exists, is not necessarily unique. As we will see the
following happens. If the covering is given by a permutation repesentation an
extension to an {\'e}tale covering of the stack will involve
an extension of the
cyclic ramification groups at the point to dihedral groups.
\end{remark}
\begin{definition}
i) A \definition {pointed curve} is a curve with only nodes as singularities
together with two disjoint finite sets of smooth points. The points
belonging to these two sets will be called \definition {cyclic}
resp.~\definition {dihedral} points.

ii) Let $(C,S_c,S_d)$ be a pointed curve, where $S_c$ is the set of cyclic and
$S_d$ the set of dihedral points. The \definition{pointed stack} associated it
is obtained from $C$ by gluing in smooth cyclic points at
$S_c$, dihedral points
at $S_d$ and cyclic points at the nodes of $C$. (Note that
the pointed curve can
be recovered from the pointed stack as the algebraic space associated to the
stack, with the sets $S_c$ and $S_d$ as the appropriate non-spatial points on
the stack.)

iii) If $(C,S_c,S_d)$ is a pointed curve then the
\definition {dual (generalised)
graph} of it is defined as follows. The vertices are the irreducible components
of $C$, for each branch of $C$ at a node one defines an oriented edge ending at
the component on which the branch lies and starting at the component on which
the other branch at the node lies (these components may of course be the
same). The edge inverse to this is the one constructed from the other branch at
the node. For each point of $S_d$ one defines an edge which equals its own
opposite starting and ending at the component on which the point lies.

iv) If $(C,S_c,S_d)$ is a pointed curve then the \definition
{dual (generalised)
graph of groups} is the graph of groups whose underlying
graph is the dual graph
of the pointed curve. The group associated to a vertex is the
fundamental group of the corresponding irreducible component of $C$ with the
union of its cyclic, dihedral and singular points removed. The group associated
to an edge corresponding to a branch at a node is the tame
local fundamental group at
the branch on the normalisation of $C$. The homomorphism into the fundamental
group at the component on which the node lies is induced by the inclusion and
the isomorphism with the tame local fundamental at the opposite branch is given
by the inverse of the canonical identification. The group associated to an edge
corresponding to a dihedral point is the tame local fundamental group, the
homomorphism into the fundamental group of the component on which it lies is
induced by the natural inclusion and the involution of the edge group given by
the fact that the edge equals its opposite is multiplication by $-1$.
\begin{remark}
To make the definition in iv) completely unambigous one would have to choose
paths and basepoints appropriately. The result, however, is independent up to
isomorphism of these choices.
\end{remark}
\end{definition}
If we want to compute the fundamental group of one of our
pointed stacks it will
of course be necessary to know when it is connected. The answer contains no
surprises.
\begin{lemma}
Let $(C,S_c,S_d)$ be a pointed curve for which $C$ is connected. Then the
associated pointed stack is connected.
\pro Let us recall that to every algebraic stack we may
construct its associated
algebraic space which is universal for maps of the stack into algebraic
spaces. Practically by construction the algebraic space associated to the
pointed stack associated to $(C,S_c,S_d)$ is $C$ itself. As the associated
algebraic space functor commutes with disjoint unions the lemma follows.
\end{lemma}
We are now ready to compute the fundamental group of a
(connected) pointed stack
associated to a pointed curve. As a point of terminology let us agree to call,
for a given abelian group $A$, the semidirect product of $A$
by a group of order
2, where the non-trivial element acts by multiplication by $-1$, the
\definition{dihedral group based on $A$}. We are also interested in the
geometric fundamental group only so we assume that hte base field is
algebraically closed.
\begin{theorem}FGR
Let $(C,S_c,S_d)$ be a pointed connected curve. Then the fundamental group of
the associated pointed stack $C''$ is isomorphic to the
fundamental group of its
dual graph.
\begin{remark}
i) I leave to the reader to write down a presentation for this group.

ii) More important than a presentation is the resulting
description of {\'e}tale
covers with structure a finite group $G$; a homorphism from the fundamental
group of each component of $C'$ to $G$ such that the
composed map from the local
fundamental group of the two branches of a singular point of
$C$ are inverses to
each other, together with a choice, for each dihedral point, of an element of
order 1 or 2 of $G$ which normalises the image of the local fundamental group
and acts by inversion on each element of that image.

iii) It is clear from the description that when $S_d$ is non-empty then the map
on fundamental groups induced from the inclusion of $C'$ in $C''$ is not
surjective. This is not altogether surprising as $C''$ is not normal.

iv) No mention of a basepoint was made in the statement of the theorem, nor of
paths from it to points close to the special points which
are needed to identify
the local fundamental groups with specific subgroups of the global one. The
reason for this is that up to isomorphism the constructions are independent of
such choices.

v) No explicit mention of which category in which the groups were to be placed
was made. The reason for this is that while we explicitly work in the algebraic
category and hence the groups are implicitly forced to be profinite, we could
have worked in the analytic category and then the theorem, interpreted in the
category of ordinary groups, would still be true.
\end{remark}
\pro Let us consider the covering of $C''$, given by $C'$
and the disjoint union
of the henselisations of all the points of $S_c$, $S_d$ as well as the
singularities of $C$ (with cyclic or dihedral points glued in). Let us assume
that we may indeed apply van Kampen's theorem to this
covering. Then by \(vK) we
get a graph $vK$ of groups whose fundamental group is the searched for
fundamental group. It then remains to show that that graph of groups has the
same fundamental group as the dual graph $D$ of groups of the pointed
curve. This can be done by considering each special point separately. A smooth
cyclic point makes $vK$ differ from $D$ by one vertex and one edge from that
vertex to the vertex representing the component on which the point lies. The
mapping from the edge group to the extra vertex group is an isomorphism so that
the vertex and edge may be removed from the graph of groups without affecting
the fundamental group. A dihedral point $d$ is responsible for one vertex and
one pair of edges between that vertex and the vertex representing the component
on which the point lies. The group at the extra vertex is the dihedral group
based on the tame fundamental group at the dihedral point and the edge group is
the tame fundamental itself. Hence this vertex and the two edges represent
(cf.~\[Se:p.~49]) an amalgamation of the dihedral group and the fundamental
group of the graph of groups with the vertex and edges removed along the tame
local fundamental group. The graph of groups $D$ differs from $vK$ by replacing
the vertex and edges associated to $d$ by a single edge which is its own
composite. Such an edge also represents the same amalgamation. Finally a node
gives rise to two vertices and four edges on $vK$ and one edge on $D$. The
verification that we again get the same contribution to the fundamental groups
is left to the reader.\note{Elaborate?}

We are left with the verification of the conditions of \(vK). This amounts
to proving the following statement.

{\it Suppose $D$ is a (not necessarily proper or connected) curve and $c$ is
either a smooth point or a node on it and we glue in a cyclic or dihedral point
to get $D'$. Let $U_1$ be the complement of $d$ (in $D'$) and let $U_2$ be the
the paullback of the stack to the Henselisation of the local ring at $d$ in
$D$. Then if we have one {\'e}tale cover over $U_1$, one
{\'e}tale cover over $U_2$ and
an isomorphism of their restrictions to the fibre product of $U_1$ and $U_2$
this set of data comes from a unique {\'e}tale cover of $D'$.}

If we take the disjoint union, $U$, of $U_1$ and $U_2$ then
general descent theory
tells that an {\'e}tale cover of this union together with descent data for the
morphism to $D'$ comes from a unique cover of $D'$. To
understand these descent data
we need to understand the fibre product of $U$ with itself over $D'$. It
consists of 3 pieces; the fibre product of $U_1$ with itself, the fibre product
of $U_2$ with itself and the fibre product of $U_1$ and $U_2$. By assumption we
have descent data on the fibre product of $U_1$ and $U_2$
and what remains to be
shown is that there is a unique way of extending these
descent data to the other
2 pieces. The fibre product of $U_1$ with itself is equal to $U_1$ as $U_1$ is
an ordinary open subset. Hence that piece takes care of itself. To understand
the fibre product of $U_2$ with itself we first consider the fibre product of
the Henselisation $\tilde R$ of the local ring at $d$ on $D$
over the local ring
$R$ at $d$. As  $\tilde R$ is Henselian this fibre product is a disjoint union
of open subsets of $\Sp {\tilde R}$. The components isomorphic to $\Sp\tilde R$
correspond to the spectrum of the tensor product of two copies of the residue
field of $\tilde R$ over the residue field of $R$, but as the latter residue
field is algebraically closed this spectrum consists of one point. Hence in the
product of $\Sp {\tilde R}$ with itself over $\Sp R$ the complement of the
diagonal is a disjoint union of copies of the spectrum of the fraction field
of $\tilde R$. Thus the mapping of this union into $D'$ factors over the fibre
product of $U_1$ and $U_2$ and thus has descent data coming from the given one,
whereas the diagonal has the obvious descent data. This proves the result
\begin{remark}
The verification of the descent condition in the corresponding analytic
situation is completely trivial, there we can let $U_2$ be a small disc with a
cyclic or dihedral point glued in and then the fibre product of it with itself
is equivalent with $U_2$ itself. The argument given in the algebraic case
simply confirms that for many practical purposes the
{\'e}tale topology behaves as
the classical one.
\end{remark}
\end{theorem}
\end{section}
\begin{section}Deformation theory.

We start by giving some general facts on equivariant
deformations of curves. Here the
assumption that the order of the group is invertible in the
base field will play
an essential role.
\begin{proposition}SD
Let $G$ be a finite group acting on a stable curve $C$ with
only ordinary double
points as singularities and assume that the order of $G$ is invertible in the
base field \k. Then the deformation problem of equivariant deformations is
formally smooth at $(C,G)$ with tangent space the
$G$-invariants on $Ext^1_{{\Cal
O}_C}(\Omega^1_C,{\Cal O}_C)$. Furthermore, the forgetful
map to the deformation
problem of non-equivariant deformations is unramified.
\pro We start by some general observations on equivariant deformations of
curves. (For the ordinary deformation theory of stable curves we refer to
\[De-Mu:\S1].) If we consider a curve $C$ with only ordinary double points as
singularities then deformations are unobstructed and the
liftings of deformations
over a small nilextension with ideal a \k-vector space $V$ is a torsor over
$Ext_C^1(\Omega^1_C,\Cal O_C)\otimes V$. All this is natural
for automorphisms of
$C$ and thus the set of liftings of deformations over a small nilextension as
above is a $G$-torsor over the $G$-module $Ext_C^1(\Omega^1_C,\Cal O_C)\otimes
V$. As multiplication on $Ext_C^1(\Omega^1_C,\Cal
O_C)\otimes V$ by the order of
$G$ is bijective all $G$-torsors are trivial and liftings always exist over
small extensions and thus over all nilimmersions. This shows that equivariant
deformations are unobstructed and the moduli stack is smooth. Furthermore, as
the tangent space at a point is a subspace of that for non-equivariant
deformations the forgetful map is unramified.
\end{proposition}
As our goal is to determine the possible equivariant
(stable) degenerations of a
smooth equivariant family of curves it makes sense to try to determine when an
equivariant stable curve has an equivariant smoothing. This is what we will do
now.
If we work in local coordinates so that an ordinary double point has ring of
functions $\k[x,y]/(xy)$ then  a miniversal deformation is given by
$\k[x,y][t]/(xy-t)$. We notice that we may choose $x$ and $y$ so that an
element of stabiliser either have them as eigenvectors or permute them up to
multiplication by a scalar. This makes the following result plausible.
\begin{proposition}ADM
Let $G$ be a finite group acting on a stable curve $C$ with
only ordinary double
points as singularities and assume that the order of $G$ is invertible in the
base field \k. There exists a deformation of the pair $(C,G)$
smoothing a singularity of $C$ iff the action of $G$ is admissible at the
point. Consequently, all singularities are equivariantly
smoothable iff they are
all admissible.
\pro
We first note that $G$ will permute the singularities and
that for each orbit of
them under $G$ that orbit is in bijection with the cosets of
$G$ with respect to
the stabiliser of a point in the orbit.  We then recall the situation in the
non-equivariant case. The tangent space to deformations of
the local singularity
is equal to the local $Ext$-group, $Ext^1_{{\Cal O}_{C,x}}(\Omega^1_C,{\Cal
O}_C)$ and the map which takes a global deformation to the
local one is just the
map from the global $Ext^1$ to the local one. As we have
seen, the tangent space
to the global equivariant problem is just the invariants under $G$ and the same
is true about but the local equivariant problem (the only difference in proof
from the global case is that there are infinitesimal automorphisms but as the
group of them has multiplication by the order of $G$ invertible no problem is
caused). More precisely, the sum of the local $Ext^1$'s in
an orbit under $G$ is
stable under $G$ and forms the induced representation of the action of a
stabiliser of one of the points of the orbit on that local $Ext^1$. The same
argument shows that equivariant deformations in the local case are unobstructed
as well. As the order of $G$ is invertible in \k, the restriction map from
global to the sum of local $Ext$-groups is surjective also on
$G$-invariants. The invariants on the sum of $Ext^1$'s over one orbit is the
invariants of the stabiliser of one point on that local $Ext^1$ and those
invariants project onto the invariants of the stabiliser of any point in that
orbit on the local $Ext^1$. This implies immediately that
the singularity $x$ is
smoothable iff $G_x$ acts trivially on the 1-dimensional space $Ext^1_{{\Cal
O}_{C,x}}(\Omega^1_C,{\Cal O}_C)$. We should thus aim to prove that a
singularity is admissible iff $G_x$ acts trivially on the local
$Ext^1$-group. The consequence in the proposition then
follows immediately. Now,
a miniversal deformation of an ordinary double point is given by
$\Spf\pow{x,y,t}/(xy-t)\to\Spf\pow t$. This shows that the tangent space to the
local deformation problem may be described invariantly as the tensor product of
the tangent spaces of the two points of the desingularisation lying over
$x$.\note{Elaborate!} This means that an automorphism of the singularity which
fixes the branches will act on the tangent space by multiplication by the
product of the character of the two branches and hence for
all those elements to
act trivially the two characters have to be inverses of each
other. On the other
hand, let us consider an automorphism \gsi{} exchanging the two branches,
choose a non-zero tangent vector of the tangent space at one of the points over
$x$ as basis for its tangent space and let its transform under \gsi{} be the
basis element for the tangent space of the other point. Having done this we see
that \gsi{} will act on the tangent space of the deformation by multiplication
by the scalar by which the square of \gsi{} acts on any of the tangent
spaces of the points lying above $x$. This means that that
square acts trivially
on the local ring at $x$ and hence in a neighbourhood of $x$.
\end{proposition}
\begin{example}
Let us consider the action of the icosahedral group, which is isomorphic to
$A_5$, on $\P^1$. There is one orbit of length 12 whose stabiliser is the Sylow
5-group of order 5. As the normaliser of that subgroup is of order 10 there is
an invariant equivalence relation on the 12 points of the orbit such that each
class contains 2 elements. Identifying each such pair gives us a curve, of
(arithmetic) genus 6, with nodes on which the icosahedral group still
acts. Each stabiliser of a node is a normaliser of a Sylow 5-group and is
dihedral so that the action of the group is admissible. Thus there exists a
smooth genus 5 curve with an action of $A_5$. Using later results of this
article we will be able to see that by taking the quotient of such a curve by
$A_5$ we get $\P^1$ and the quotient map is ramified at 4 points on $\P^1$. It
is also easy to figure out the ramification behaviour for this map. In
this way we get an algebraic proof for the existence of a particular ramified
covering of $\P^1$ (whose existence, of course is ensured by
Riemann's existence
theorem).
\begin{remark}
It should be emphasized that this is not a viable strategy for giving an
algebraic proof of Riemann's existence theorem in general; the irreducible
components of a degeneration will in general be more complicated than just a
rational curve.
\end{remark}
\end{example}
We will also need to know what happens with the action of a stabiliser under
smoothing. The result is given in the following lemma.
\begin{lemma}FDEF
Let the finite group $G$ act on a curve $C$. Let $c\in C$ be an admissible
singularity, assume that the stabiliser of $c$ acts faithfully on $c$ and let
$\Cal C\to B$ be an equivariant deformation of $C$ with $b\in B$ corresponding
to $C$. We will say that a geometric fibre $C'$ of $\Cal
C\to B$ over a generisation
$b'\in B$ of $b$ for which $c$ is not the specialisation of a singularity of
$C'$ is a \definition {smoothing} of $C$. A point on $C'$ of a having $c$ as a
specialisation will be referred to as a point specialising to $c$.

i) If $c$ is a cyclic node then there will be no point with a non-trivial
stabiliser specialising to $c$.

ii) If $c$ is a dihedral node then there will be 4 or 2 orbits of points on a
smoothing specialising to $c$ whose stabilisers are dihedral involutions,
depending on whether the order of the stabiliser is divisible by 4 or not. No
other points with non-trivial stabilisers specialise to $c$.
\pro We first want to show that the statement is local around $c$. Indeed, for
i) this is obvious and for ii) one needs to prove that if two neighbouring
fixpoints are conjugate they are conjugate under the action
of the stabiliser of
$c$. This however is clear as an element taking one fixpoint to another has to
take the irreducible component of the fixpoint locus at the first point to the
irreducible component of the fixpoint locus at the other
point. As both of these
components pass through $c$ the element has to be in the
stabiliser. Thus we may
localise around $c$ and we may also pass to a miniversal deformation. Hence we
are reduced to looking at an admissible action on
$\pow{x,y,t}/(xy-t)$ over $\pow t$
and we may also assume that an element of $G$ either acts by multiplying $x$ by
a scalar and by multiplying $y$ by the inverse of that scalar or permutes acts
by $x\mapsto \gla y$ and $y\mapsto \gla^{-1}x$ for some non-zero \gla. If $g\in
G$ acts on $x$ by multiplication by $\gla\ne1$, then the fixpoint locus of $g$
is defined by the ideal generated by $gx-x=(\gla-1)x$ and
$gy-y-(\gla^{-1}-1)y$. As \gla{} is different from 1 this equals the ideal
generated by $x$ and $y$ and thus the fixpoint locus equals $c$. This proves
i). If $g$ instead takes $x$ to $\gla y$ and $y$ to $\gla^{-1}x$ then the
fixpoint locus is defined by $\gla y-x$ and $\gla^{-1}x-y$ i.e.~by $\gla
y-x$. Dividing out $\pow{x,y,t}/(xy-t)$ by this relation gives
$\pow{x,y,t}/(\gla y^2-t)$ which is finite flat of degree 2 over $\pow t$ and
thus gives  2 fixpoints for neighbouring values of $t$. The
dihedral involutions
form 2 or 1 conjugacy classes depending on whether the order of the dihedral
group is divisible by 4 or not. This gives 4 resp.~2 orbits of fixpoints.
\end{lemma}
\end{section}
\begin{section}On the boundary of the Hurwitz scheme

In this section we will start tying things together. We
begin with a result which
shows that our constructions do indeed capture the whole boundary of smooth
curves with an action of a finite group.
\begin{definition}
For a finite group $G$ and an integer $g>1$ we denote by $\Cal M_g(G)$ the
moduli stack of smooth, connected curves of genus $g$
together with an action of $G$. We
also denote by $\overline{\Cal M}_g(G)$ the moduli stack of
stable, connected curves of
(arithmetic) genus $g$ together with an admissible action of $G$ on the curve.
\end{definition}
\begin{remark}
In general, $\Cal M_g(G)$ is of course not connected even if one would restrict
to only actions which are faithful (which would be a reasonable thing to do).
\end{remark}
We then get the following result.
\begin{theorem}
$\overline{\Cal M}_g(G)$ is a smooth and proper algebraic stack with $\Cal
M_g(G)$ as a dense open substack and the forgetful map
$\overline{\Cal M}_g(G)\to
\overline{\Cal M}_g$ is unramified.
\pro The fact that this stack is formally smooth follows
directly from \(SD) and
that $\Cal M_g(G)$ is a dense open subset follows from
\(ADM). That the stack is
of finite type follows directly from the fact that $\Cal M_g$ is and the
finiteness of the automorphism group of a stable curve of genus $>1$. To show
properness we use the valuative criterion. Hence we assume
that we have a stable
curve with an action by $G$ over the fraction field $K$ of a discrete valuation
ring $R$. By the stable reduction theorem after possibly extending $K$ we may
assume that there is a unique extension of the curve to a flat family over $R$
with stable fibres. By the uniqueness the action of $G$ will extend to this
family and as the generic fibre is dense in the total space this extension is
unique. This shows properness. The fact that the forgetful map is unramified
follows from \(SD).
\end{theorem}
Having the results of the two previous sections in mind it should be clear that
the data given by an admissible action of a group $G$ on a
curve with only nodes as
singularities is very closely related to the data given by a pointed curve and
{\'e}tale covers of the associated stack. The precise relation is as follows.
\begin{theorem}
Let the group $G$ act admissibly on the curve $C$ with only nodes as
singularities. Let $C'$ be the algebraic stack obtained from $C$ by gluing in
cyclic nodal points at all the nodes and smooth cyclic points at all smooth
points of $C$ having a nontrivial stabiliser group under the action of
$G$. The action of $G$ on $C$ extends (uniquely) to an action of $G$ on
$C'$. Let $C''$ be the stack quotient of $C'$ by $G$ (so that the quotient map
$C'\to C''$ is {\'e}tale). Let $C_1$ be the quotient of $C$
under $G$ (as an ordinary
scheme) and let $C_1''$ be the stack obtained from $C_1$ by gluing in smooth
cyclic points at all points of $C_1$ below a smooth point with non-trivial
$G$-stabiliser, dihedral points at smooth points below a node on $C$ and nodal
cyclic points at nodes of $C_1$. Then then natural
isomorphism between $C''$ and
$C_1''$ with their non-scheme points removed extends (uniquely) to an
isomorphism between $C''$ and $C_1''$. In particular, $C$ is the scheme
associated to the {\'e}tale $G$-cover $C'$ of $C''_1$.
\pro
As the morphism is defined outside of the special points only the existence of
an extension needs to be proved and that is a question local around the
singularities of $C$ and the points of $C$ with non-trivial
$G$-stabilisers. From this one immediately reduces down to something local at
points of $C$ under the action of the $G$-stabiliser of the point. Hence we may
consider the strict Henselisation $R$ of a point and an admissible action of a
group on that point. In all cases we have to compare the result obtained by
first gluing in a non-scheme point and then taking the stack quotient by the
group and that obtained by first taking the scheme quotient
and then gluing in a
non-scheme point. Let us first consider the case of a smooth point with a
(necessarily) cyclic stabiliser. Gluing in a cyclic point means passing to the
normalisation $\tilde R$ of $R$ in the maximal tamely ramified extension and
taking the stack quotient by the Galois group of that extension. Then taking a
further stack quotient by a cyclic group $G$ is clearly the same thing as
considering $\tilde R$ the normalisation of the maximal tamely ramified
extension of the invariant ring of $G$ and then taking the stack quotient. On
the other this is exactly taking the invariants and then taking the stack
quotient by the Galois group of the maximal tamely ramified extension acting on
the normalisation of these invariants. This shows the result
for the action on a
smooth point. The case of the action of a cyclic group on a node is completely
similar. The case of a dihedral group acting on the node is slightly
different. To get that case one needs to see that the following two
constructions give the same result. 1'st construction: Start with a strictly
Henselian local ring with an ordinary node as closed point. Then consider the
subring of the normalisation in the product of the maximal tamely ramified
extensions of the two fraction fields consisting of pairs $(r,s)$ congruent to
each other modulo the respective maximal ideals (as in the construction of the
gluing in of nodal cyclic point above). 2'nd construction:
Consider an action of
a dihedral group acting admissibly on the Henselisation of a node and take the
invariant ring which is a regular ring. Then take the
normalisation of this ring
in the maximal tamely ramified extension and consider the
subring of the product
of this ring with itself consisting of the pairs $(r,s)$ which are equal modulo
the respective maximal ideals (as in the construction of the gluing in of a
dihedral point). Now, we may first take invariants under the normal cyclic
subgroup to reduce to the case of a dihedral group of order 2. Then it is clear
that the normalisations of the constructed rings are the same and hence so are
the rings themselves being defined as subrings of their normalisations by
identical conditions.
\end{theorem}
An immediate consequence of this result is the principal technical
result of this paper. To formulate it let us agree to the following definition.
\begin{definition}
Let a finite group $G$ act admissibly on a curve $C$ with only nodes as
singularities. We give the quotient curve $C/G$ the
structure of a pointed curve
by letting the cyclic points be the smooth points below smooth points with
non-trivial $G$-stabiliser and the dihedral points lying below nodal points. We
will call this pointed curve the \definition{pointed curve
quotient} of $C$ by $G$.
\end{definition}
This definition then, naturally, leads to the following
corollary of the theorem.
\begin{corollary}
Let $(C,S_c,S_d)$ be a pointed curve with $S_c$ the set of cyclic points and
$S_d$ the set of dihedral points. For a finite group $G$ there is a natural 1-1
correspondence between admissible actions of $G$ on a curve with only nodes as
singularities together with an isomorphism of the pointed quotient with
$(C,S_c,S_d)$ and {\'e}tale $G$-covers of the pointed stack associated to
$(C,S_c,S_d)$. In particular such actions correspond to conjugacy classes of
homomorphisms of the fundamental group of the pointed stck to $G$.
\pro That one gets an {\'e}tale cover of the pointed stack
is the content of the
theorem. To go from an {\'e}tale cover of the stack to a
curve one simply looks at
the algebraic space associated to the stack cover.
\end{corollary}
\begin{remark}
The explicit consequence of this corollary combined with \(FGR), namely a
concrete description of actions of $G$ on curves with only nodes as
singularities in terms of the pointed curve quotient, can be formulated without
ever mentioning algebraic stacks. The reason for using algebraic stacks is that
a lot of tedious verifications are avoided by referring instead to the general
fact that the category of {\'e}tale coverings of a connected
stack form a Galois
category. Let me mention however in concrete terms how one constructs the curve
cover associated to a representation of the fundamental
\gpi\ group of the generalised
graph of groups associated to a pointed curve $(C,S_c,S_d)$ in a group
$G$. First because the fundamental group of each irreducible component of $C$
minus the union of the singularities of $C$, $S_c$ and $S_d$ maps into \gpi, we
get a ramified $G$-cover of such a component. The
amalgamation condition for the
two branches of a node of $C$ implies that there exists a $G$-equivariant
identification of the points of the respective ramified covers over the two
branches. One may thus equivariantly pair these points off
to from nodes mapping
down to the appropriate node of $C$. Similarly, the dihedral
involution given by
the amalgamation at a dihedral point gives a pairing off of points above a
dihedral point to construct nodes above it.
\end{remark}
If we now go back to the motivating example of the Hurwitz
scheme of ramified covers
of the projective line and their descriptions in terms of such covers which are
Galois we see that we have indeed obtained a description of possible
degenerations of such covers: We start with a family of covers of the line, we
form the family of Galois covers and then (after possibly extending the base)
extend it to a complete family of stable curves with the
Galois group $G$ acting
on it. We then finally take the quotient of this family by the stabiliser of a
point in the permutation representation of $G$ given by the original family
which gives us a compactification of the original family. Compactifications of
Hurwitz families have, of course, been considered previously (for instance in
\[Ha-Mu]) and then they are done directly in terms of the covering of the
projective line. We will now make a few comments on the comparisons of the
approach given in this article and these more direct approaches.
\begin{remark}
In the preceding paragraph, as well as in the paragraphs to follow the precise
relations between moduli problems are glossed over. The aim of the present
article is to give a convenient description of the topological types of
degenerations rather than to describe a suitable moduli problem. Thus we are
primarily interested in individual curves rather than families of curves. The
precise relations between moduli problems are therefore left to the interested
reader. (Some of the differences are expounded upon in what follows.)
\end{remark}
Let us now go on to make a rough comparison of the present
approach with that of
\[Ha-Mu]. Thus consider a stable curve $C$ with a (faithful) action of the
finite group $G$ such that the quotient of $C$ by $G$ is of genus 0. We may
consider the fixpoint loci of non-trivial elements of $G$ and then throw away
the parts corresponding the branch preserving actions on singular points. We
then consider the quotient of these loci as a subscheme on the quotient of $C$
by $G$. By construction the union of them is a subscheme with support in the
smooth part and it follows from the calculation done in the proof of \(FDEF)
complemented by a similar calculation  for a cyclic point that this sbscheme is
of length 2 at a dihedral point and of length 1 at a smooth.
By the calculations
of \(FDEF) and a similar calculation at a smooth cyclic point this construction
behaves well in families so that the quotient of $C$ by $G$ naturally comes
equipped with the structure of a relative divisor with support in its smooth
part and with local multiplicity at most 2. On the other hand in the setup of
\[Ha-Mu] a genus 0 curve with a relative divisor supported in the smooth part
which has everywhere local multiplicity 1 is used as a starting point. The
general relation between these two structures would need some elaboration but
let me explain it in a ``miniversal'' situation. We
therefore assume that we have a
1-dimensional family of stable curves with $G$-action, a point $c$ on a fibre
$D$ which is a dihedral node on that fibre (and for simplicity assume it to be
the only one) and such that the divisor of fixed points on the quotient of the
family by $G$ has two irreducible components at the image $d$ of $c$ both of
them intersecting the fibre to which $d$ belongs as well as the other component
transversally. (This situation can be obtained by a
quadratic base change from a
miniversal deformation, which follows from the calculations of the proof of
\(FDEF).) We then blow up the point $d$. This adds a projective line to $C$,
giving the new fibre $C'$, meeting one other component in a smooth point. The
divisor of fixpoints will meet this new component in 2 points both of which lie
only on that component. Thus $C'$ is a genus 0 curve with a relative divisor
supported in the smooth part and which has everywhere local multiplicity 1. It
is also a pointed curve whose cyclic points are the cyclic points of $C$ plus
the 2 points on the new component and with no dihedral points. Its fundamental
group is obtained from that, $\Pi$ say, of the pointed curve which is $C$ but
with $d$ as a cyclic point by amalgamating the local fundamental group at $d$
with the local fundamental group at one point of $\P^1$ minus 3 points. The
fundamental group of $C$ is the amalgamation of $\Pi$ with a dihedral group
based on the local fundamental group at $d$ along that local fundamental
group. There is an obvious group homomorphism from the fundamental group of
$\P^1$ minus 3 points to the dihedral group based on one of its local
fundamental groups taking generators of the two other local fundamental groups
to dihedral involutions. This morphism glues to a surjective map from the
fundamental group $\Pi_1$ of the pointed curve $C'$ to that
of the pointed curve
$C$, $\Pi_2$ say. Now if think of the action of $G$ on $D$ as a group
homomorphism from $\Pi_2$ to $G$ we may compose that with the surjection from
$\Pi_1$ to $\Pi_2$ and hence get a $G$-covering $D'$ of $C'$. The curve $D'$
maps equivariantly to $D$ and may be thought of as being obtained from $D$ by
replacing each dihedral node by a projective line which meets the two branches
of the node in 0 and $\infty$ and for which the stabiliser in $G$ of the node
acts in the standard way on the projective line (with the dihedral involutions
exchanging 0 and $\infty$). In the situation of Harris and Mumford we also have
a transitive permutation representation of $G$ and by dividing $D$ by the
stabiliser of a point in that representation we get a curve $E$ which is an
admissible covering of $C'$ (except that to get only simple
ramification of this
covering one needs some conditions on the permutation representation, however
Harris and Mumfords construction immediately generalises to the case when this
is not assumed). Thus very roughly speaking the Harris-Mumford construction
avoids dihedral stabilisers by blowing up those points (all the stabilisers on
$D'$ are cyclic).
\end{section}
\begin{section}Codimension 1 and some examples

We will now discuss the simplest possible degenerations, those of codimension 1
in the moduli space and look at a few examples. We saw in our study of the
deformation theory of admissible actions that any
equivariant deformation of the
singularities of a curve could be realised by a global equivariant
deformation. This means that we can always equivariantly smooth one orbit of
singularities while keeping the singularities in other orbits. This means that
one may smooth all but one orbit and hence the points of codimension 1 at the
boundary of moduli space will correspond to a single orbit of
singularities. Otherwise put, the pointed quotient curve
will have either 1 node
or 1 dihedral point and then the rest of the special points will be cyclic. To
simiplify matters let us assume that we are in the case motivating us, namely
that the quotient curve is of genus 0. Then it will either be just the
projective line with 1 dihedral and a certain number, $n$
say, (at least 2 to make it
stable) of cyclic points, or it will have as irreducible
components 2 projective
lines meeting in a node and then a certain number, $n_i$,
$i=1,2$ say, of cyclic
points on each component (at least 2 on each to make the curve stable). In the
first case the fundamental group will have generators $\gsi,x_0,x_1,\ldots,x_n$
with relations $\gsi^2=e$, $\gsi x_0\gsi=x^{-1}_0$ and $x_0x_1\cdots x_n=e$ and
in the second case generators $x_1,x_2,\ldots,x_{n_1}$ and
$y_1,y_2,\ldots,y_{n_2}$ with relations $x_1x_2\cdots x_{n_1}=e$, $y_1y_2\ldots
y_{n_2}=e$ and $x_{n_1}y_1=e$. An equivalent, but sometimes more convenient,
presentation in the dihedral case is obtained by setting $\tau:=\gsi x_0$ and
then we get generators $\gsi,\tau,x_1,\ldots,x_n$ with relations
$\gsi^2=\tau^2=e$ and $\gsi\tau x_1\cdots x_n=e$. In the dihedral case, we will
have use for the following observation. If we consider the normalisation of the
pointed stack associated to the pointed curve consisting of the line with the
chosen dihedral point and $n$ chosen cyclic points then the map on fundamental
groups induced by the normalisation map is just the inclusion of the group
generated by $x_0,x_1,\ldots,x_n$ into the one generated by
$\gsi,x_0,x_1,\ldots,x_n$, Thus if we have a finite group $G$ and elements
$s,g_0,g_2,\cdots g_n$ with $s^2=e$, $sg_0s=g_0^{-1}$ and $g_0g_1\cdots g_n=e$
and let $C$ be the corresponding curve with $G$-action whose quotient is the
projective line then the normalisation of $C$ as a (ramified) $G$-covering of
the line is classified by the elements $g_0,g_1,\cdots,g_n$ thought of as
describing a $G$-covering of the line unramified outside of
the given dihedral and
cyclic points.

Let us now consider our previous example of the icosahedral group acting on the
projective line in this light. Thus we pick elements $g_0$, $g_1$ and $g_2$ of
order 5, 2 and 3 respectively whose product is one in $A_5$. This gives a
$A_5$-covering of $\P^1$ ramified at 0, 1 and $\infty$. Now
the normaliser of a Sylow
5-subgroup of $A_5$ is dihedral so we may pick elements $s$ and $t$ of it of
order 2 s.t.~$st=g_0$. Now consider the pointed curve which is $\P^1$ with 0 as
a dihedral point and 1 and $\infty$ as cyclic points. The collection
$s,t,g_1,g_2$ gives  an {\'e}tale $A_5$-cover of the associated stack whose
normalisation is the $A_5$-cover constructed using $g_0,g_1,g_2$. This cover is
the one considered in the previous example. Now, when smoothing it we get a
covering of $\P^1$ ramified in 4 points corresponding to the tuple
$s,t,g_1,g_2$. Using our results on the d{\'e}vissage of the
action of the group on
cohomology we get a quick way of computing the character of the action of $A_5$
on these smooth genus 6 curves; it is simply 2 times the non-trivial character
of degree 1 of the subgroup of $A_5$ of order 10 induced up to $A_5$. (This
character can of course be computed also from the ramification behaviour and
Lefschetz fixpoint formula). The same example can also be obtained by dividing
$s,t,g_1,g_2$ in the groups $s,t,(st)^{-1}=g_0^{-1}$ and $st=g_0,g_1,g_2$ which
gives a cover of two crossing $\P^1$'s with two cyclic
points on each component,
the covering of one component is the icosahedral, the covering of the other is
given by the action of a dihedral group of order 10 on $\P^1$.

A similar example is obtained by choosing elements of order 7, 2, and
3 in $PSL(2,\F_7)$ whose product is the identity and then use that the
normaliser of the Sylow 7-group is dihedral of order 14. Our
starting curve -- a
$PSL(2,\F_7)$-cover of $\P^1$ ramified at 3 points with ramification groups of
order 7, 2, and 3, is then the Klein curve of genus 3. In the end we get a
$PSL(2,\F_7)$-covering of genus 15 of $\P^1$ ramified in 4 points with
ramification groups of order 2, 2, 2, and 3.
\end{section}
\begin{bibliography}
\[De-Mu]:\by P. Deligne, D. Mumford\paper On the irreducibility of the moduli
space of curves \jour Publ. I.H.E.S.\vol 36\yr 1969\pages 75 -- 110

\[Ha-Mu]:\by J. Harris, D. Mumford\paper On the Kodaira dimension of the moduli
space of curves\jour Inv. math.\yr 1982\vol 67\pages23 -- 86

\[Ka-Ma]:\by N. Katz, B. Mazur\book Arithmetic moduli of elliptic curves\publ
Princeton Univ. Press \publadr Princeton\yr 1985

\[Se]:\by J.-P. Serre\book Trees \publ Springer Verlag\publadr Berlin\yr 1980

\end{bibliography}
\end{document}